\documentclass{amsart}

\usepackage{graphicx}
\usepackage{authblk}
\usepackage{upgreek}
\usepackage{hyperref}
\usepackage{float}

\newtheorem{theorem}{Theorem}[section]
\newtheorem{prop}[theorem]{Proposition}

\newtheorem{lemma}[theorem]{Lemma}

\DeclareMathOperator{\dist}{dist}
\DeclareMathOperator{\ff}{\mathfrak{f}}
\DeclareMathOperator{\fa}{\mathfrak{a}}
\DeclareMathOperator{\fu}{\mathfrak{u}}
\DeclareMathOperator{\fp}{\mathfrak{p}}
\DeclareMathOperator{\fq}{\mathfrak{q}}
\DeclareMathOperator{\fg}{\mathfrak{g}}
\DeclareMathOperator{\fk}{\mathfrak{k}}

\DeclareMathOperator{\lm}{\mathcal{L}}

\DeclareMathOperator{\real}{\mathbb{R}}
\DeclareMathOperator{\r2}{\mathbb{R}^2}
\DeclareMathOperator{\integer}{\mathbb{Z}}

\usepackage{xcolor}

\usepackage[foot]{amsaddr}
\usepackage{appendix}
\begin{document}

\title[Examples of singularity models for $\mathbb{Z}$/2 harmonic 1-forms and spinors]{Examples of singularity models for $\mathbb{Z}$/2 harmonic 1-forms and spinors in dimension three}

\author{C. H. Taubes$^\dag$}
\author{Y. Wu$^\lozenge$}

\address{$^\dag$C. H. Taubes was supported in part by the NSF (DMS 1708310) and by the Stanford University Henri Poincar\'e Distinguished Visiting Professor Fund.  $^{\dag\lozenge}$Both authors thank the Department of Mathematics at Stanford for its hospitality.}

\begin{abstract}
We use the symmetries of the tetrahedron, octahedron and icosahedron to construct local models for a $\mathbb{Z}/2$ harmonic 1-form or spinor in 3-dimensions near a singular point in its zero loci.  The local models are $\mathbb{Z}/2$ harmonic 1-forms or spinors on $\mathbb{R}^3$ that are homogeneous with respect to the rescaling of $\mathbb{R}^3$ with their zero loci consisting of 4 or more rays from the origin. The rays point from the origin to the vertices of a centered tetrahedron in one example, and to those of a centered octahedron and a centered icosahedron in two others. 
\end{abstract}

\maketitle

\section{Introduction}
	Suppose in what follows that $X$ denotes a smooth, oriented, Riemannian 3-manifold.  A $\mathbb{Z}/2$ harmonic 1-form on $X$ consists of a data set $(Z,\mathcal{I}, v)$ whose constituents are as follows:  What is denoted by $Z$ signifies a closed subset of $X$ with Hausdorff dimension at most 1.  What is denoted by $\mathcal{I}$ signifies an associated $\mathbb{R}$ bundle to a principle $\mathbb{Z}/2$ bundle over $X - Z$, hence a real line bundle over this domain.  What is denoted by $v$ signifies a closed and coclosed 1-form with values in $\mathcal{I}$ whose norm extends over $Z$ to define a H\"older continuous function on $X$ that vanishes on $Z$.  To say that $v$ is closed and coclosed is to say that it obeys the equations
\begin{eqnarray*}\label{eqn:1.1}
dv = 0  \  \text{ and } \  d * v = 0
\end{eqnarray*}
\hfill (1.1)\\
with $*$ denoting the metric’s Hodge dual operator.  (With regards to taking derivatives of sections of $\mathcal{I}$:  Derivatives are defined over any given ball in $X-Z$ by choosing an isometry to identify $\mathcal{I}$ over the ball with the product $\mathbb{R}$-bundle.)  

	A $\mathbb{Z}/2$ harmonic spinor over $X$ consists of a data set $(Z,\mathcal{I}, s)$ with $Z$ and $\mathcal{I}$ as before and with $s$ being an $\mathcal{I}$-valued spinor on $X$ (with respect to a chosen spin structure) that obeys the Dirac equation on $X-Z$ and whose norm extends over $Z$ as a H\"older continuous function on $X$ that vanishes on $Z$.		
		
	These $\mathbb{Z}/2$ harmonic gadgets (1-forms and spinors) are of interest because they characterize in part the behavior of non-convergent sequences of solutions to certain first-order gauge theory equations:  The $\mathbb{Z}/2$ harmonic 1-forms characterize (in part) the behavior of non-convergent sequences of equivalence classes of flat $Sl(2;\mathbb{C}$) connections on $X$ (see \cite{taubes2013psl,taubes2015corrigendum}); and $\mathbb{Z}/2$ harmonic spinors characterize in part the behavior of non-convergent sequences of equivalence classes of solutions to the 2-spinor generalization of the Seiberg-Witten equations (see \cite{haydys2015compactness}).  	
	
As explained by Takahashi (see \cite{takahashi2015moduli} and \cite{takahashi2017index}) and elaborated on by Donaldson \cite{Donaldson}, there is a well behaved moduli space of $\mathbb{Z}/2$ harmonic 1-forms and spinors near any given $(Z,\mathcal{I}, v$ or $s$) in the case when $Z$ is a $\mathcal{C}^1$ embedded submanifold in $X$.  But, it is not known a priori that this is always the case.  Even so, a theorem of Zhang \cite{zhang2017rectifiability} says that $Z$ is always rectifiable and that it always has finite 1-dimensional Hausdorff measure. Thus, it has a dense subset with the structure of a $\mathcal{C}^1$ submanifold (see also \cite{taubes2014zero}).  

Supposing that $Z$ is not everywhere a $\mathcal{C}^1$ submanifold, then there are local models for its singular points, which are $\mathbb{Z}/2$ harmonic 1-forms (or spinors) on $\mathbb{R}^3$ that are homogeneous with respect to coordinate rescalings.  To elaborate:  A coordinate rescaling is a linear diffeomorphism of $\mathbb{R}^3$ (with $\mathbb{R}^3$ viewed as a vector space) that sends any given vector (call it $x$) to $\lambda x$ with $\lambda$ being a positive number.  A $\mathbb{Z}/2$ harmonic 1-form or spinor on $\mathbb{R}^3$ is homogeneous with respect to coordinate rescalings when the conditions listed below in (1.2) are met.  (In the case of spinor, the pull-back is defined via a suitable lift to the spin bundle of the action of the group of rescaling diffeomorphism.)

\begin{itemize}\it
\item	The set Z is a finite union of rays from the origin and thus mapped to itself by any coordinate rescaling diffeomorphism.  
\item The pull-back of $\mathcal{I}$ via any coordinate rescaling diffeomorphism is isomorphic to $\mathcal{I}$.  
\item The pull-back of the 1-form $v$ or spinor $s$ by the rescaling defined by any given positive number $\lambda$ has the form $\lambda^\alpha v$ or $\lambda^\alpha s$ with $\alpha$ being independent of $\lambda$.
\end{itemize}
\hfill\rm (1.2)

The simplest example of homogeneous $\mathbb{Z}/2$ harmonic 1-form follows:  Let $x =  (x_1, x_2, x_3)$ denote Euclidean coordinates for $\mathbb{R}^3$ and let $z$ denote the complex coordinate $z = x_1 + i x_2$.  Set $v$ to be the real part of $\sqrt{z}dz$; or the real part of $z^k\sqrt{z}dz$ with $k$ being a positive integer.  The set $Z$ in these cases is the $x_3$-axis.  This example supplies the local model for the non-singular part of the vanishing loci of a $\mathbb{Z}/2$ harmonic 1-form on a Riemannian 3-manifold.   There is a similar local model for the non-singular part of vanishing loci of a $\mathbb{Z}/2$ harmonic spinor on a Riemannian 3-manifold where the spinor has the form $z^{k+{1\over2}} s$ with $k$ being a non-negative integer and $s$ being a suitable constant spinor.

This article supplies a handful of local models for singular loci of $\mathbb{Z}/2$ harmonic 1-forms and spinors.  By way of a look ahead, the versions of $Z$ for these models comprise 4 or more rays from the origin.  The simplest case has $Z$ being 4 rays, the rays from the origin through the vertices on the $|x| = 1$ sphere of an inscribed regular tetrahedron.  Another example has $Z$ being the rays from the origin through the vertices on the $|x| = 1$ sphere of an inscribed, regular icosahedron.  Yet another example has $Z$ being the 20 rays from the origin through the midpoint of the faces of this same icosahedron.   To the authors’ knowledge, these are the first examples of homogeneous $\mathbb{Z}/2$ harmonic 1-forms (and spinors) that are not $SO(3)$ rotations of the ones that are described in the preceding paragraph. The appendix to this article proves a proposition to the effect that the only homogeneous, $\mathbb{Z}/2$ harmonic 1-forms or spinors on $\mathbb{R}^3$ with the set $Z$ being a union of just two rays from the origin are those where the rays are antipodal.  (This is the case for those in the preceding paragraph.)

The example given here of a homogeneous, $\mathbb{Z}/2$ harmonic 1-form on $\mathbb{R}^3$ where $Z$ is the union of 4 rays from the origin with no two being pairwise colinear can be viewed as an $\mathbb{R}$-invariant, homogeneous, $\mathbb{Z}/2$ harmonic 1-form on $\mathbb{R}^4$ by viewing $\mathbb{R}^4$ as $\mathbb{R}^3 \times \mathbb{R}$.  Viewed in this light, the example is one where the 4-dimensional analog of $Z$ is the union of 4 half-planes in $\mathbb{R}^4$ that share a common boundary line but with no two being coplanar.  This geometry for the 4-dimensional version of $Z$ was listed in \cite{taubes2014zero} as one of the allowed (in principle) forms for the versions of $Z$ that can appear in a 4-dimensional, $\mathbb{Z}/2$ harmonic 1-form singularity model.  Examples in $\mathbb{R}^4$ with $Z$ being a union of multiple (full) planes through the origin in $\mathbb{R}^4$ are explicitly described in \cite{taubes2014zero}.

Our singularity models are constructed using solutions to a second order differential equation on the $|x| = 1$ sphere.  To say more about this related problem:  Let $N$ denote an even, positive integer and let $Z\equiv \{p_1,\ldots, p_N\}$ denote a set of $N$ distinct points in $|x| =1$ sphere (this sphere is denoted by $S^2$).   
The fundamental group of $S^2 - Z$ is a free group of rank $N-1$.  As such, it is convenient to fix a set of $N$ generators (the set is denoted by $\{\gamma_1, \ldots, \gamma_N\}$ for the abelianization subject to the one relation $\gamma_1 \cdots \gamma_N = 1$ with any given $\gamma_{k}$ represented by a small radius circle about the corresponding point $p_k$ whose interior contains only $p_k$.  
Let $\iota$ denote the homomorphism from $\pi_1(S^2 - Z)$ to $\mathbb{Z}/2$ that sends each $\gamma_k$ to $-1$.  This homomorphism defines a principle $\mathbb{Z}/2$ bundle over $S^2 - Z$.  (This principle bundle is the complement of the branch points in the 2-sheeted branched cover of $S^2$ with branching loci $Z$.) 
 We use $\mathcal{I}$ in what follows to denote the real line bundle over $S^2 - Z$ that is associated to this same $\mathbb{Z}/2$ principle bundle.  The points in $Z$ are said to be the points of discontinuity of $\mathcal{I}$.  

A section $\ff$ of $\mathcal{I}$ over ($S^2-Z$) can be viewed as a function on $S^2-Z$, which is defined at any given point up to multiplication by $\pm 1$.  The sign changes from $+1$ to $-1$ when circling any given generating loop $\gamma_k$.  Locally on $S^2-Z$, a section of $\mathcal{I}$ can be viewed as just an ordinary function.  Because of this, the exterior derivative of functions on $S^2$ (which acts locally by taking first derivatives) gives a map from sections of $\mathcal{I}$ over $S^2-Z$ to $\mathcal{I}$ valued 1-forms.  The exterior derivative of a section $\ff$ is denoted by $d\ff$.  (Keep in mind that $d\ff$ is not defined at the points in $Z$.)  Likewise, the Laplacian on functions on $S^2$ (which acts locally by taking second derivatives) sends any given sections of $\mathcal{I}$ over the domain $S^2 - Z$ to sections of $\mathcal{I}$.  This is denoted by $\Delta$.  A section $\ff$ is said to be an eigensection of $\Delta$ when $\Delta \ff = -\mathcal{E}\ff$ with $\mathcal{E}$ being a real number. 

Of principle interest here are eigensections with the following property:\\

\begin{center}\it
The norms of $\ff$ and $d\ff$ extend over $Z$ to define H\"older continuous functions on $S^2$ that vanishes at the points in $Z$.
\end{center}
\rm\hfill(1.3)

By way of an example:  Let $(x_1, x_2, x_3)$ again denote the Cartesian coordinates on $\mathbb{R}^3$.  Take $Z$ to be the set $\{(0, 0, 1), (0, 0, -1)\}$.  Let $a$ denote a non-zero complex number, let $k$ denote a positive integer and set $\ff$ to be the restriction to $S^2-Z$
of the real part of $a(x_1+ix_2)^{k+{1\over2}}$.
   
Eigensections that obey (1.3) will be used to construct the local models for the singularities of $\mathbb{Z}/2$ harmonic 1-forms and spinors on 3-manifolds.  To this end, let $Z \equiv \{p_1, \ldots, p_{2N}\}$ denote a set of $2N$ points on the $|x| = 1$ sphere in $\mathbb{R}^3$ and also the rays from the origin through these same points.  Meanwhile, let $\mathcal{I}$ denote the real line bundle over $S^2-Z$ of the sort that is described above and also its pullback to the complement in $\mathbb{R}^3$ of the eponymous set of rays via the map (denoted by $\pi$) that sends any given point $x$ to ${x\over |x|}$.  Set $\ff$ to denote a section of $\mathcal{I}$ defined over $S^2-Z$ that obeys $\Delta\ff = -\mathcal{E}\ff$ and also the conditions in (1.3).  Now let $v$ denote the 1-form $d(|x|^\alpha\pi^*\ff)$ with $\alpha$ denoting  ${1\over 2}(1+(1+4\mathcal{E})^{{1\over2}})$. The data consisting of the ray set $Z$ in $\mathbb{R}^3$, the line bundle $\mathcal{I}$ (pulled back from its eponymous line bundle on $S^2 - Z$ via $\pi$) and $v$ defines a homogeneous $\mathbb{Z}/2$ harmonic 1-form data obeying (1.2).  

 Conversely, any such $\mathbb{Z}/2$ harmonic 1-form data set on $\mathbb{R}^3$ that is described by (1.2) has the form $v = d(|x|^\alpha \pi^*{\ff})$ with $\ff$ being an eigensection of the $S^2$ Laplace operator acting on sections of the restriction of the line bundle from the second bullet in (1.2) to the complement in $S^2$ of the rays that comprise the set $Z$ in the top bullet of (1.2). 
 
	Local models for the $\mathbb{Z}/2$ harmonic spinor singularities can also be obtained from data $(Z, \mathcal{I},\ff)$ on $S^2$:  Let $\mathbb{D}$ denote the Dirac operator acting on sections of the product spinor bundle over $\mathbb{R}^3$, the bundle $\mathbb{R}^2 \times \mathbb{C}^2$. When written using Cartesian coordinates, $\mathbb{D}$ is the $2 \times 2$ matrix operator
$$\mathbb{D} = \begin{pmatrix} i{\partial \over \partial x_3} & i {\partial \over \partial x_1} + {\partial \over \partial x_2}\\
i{\partial \over \partial x_1} - {\partial \over \partial x_2} & -i{\partial \over \partial x_3}\end{pmatrix}.$$
\hfill (1.4)\\
Fix a non-zero, constant spinor which will be denoted by $s_0$.  The data consisting of the rays from the origin through the points of $Z$, the $\pi$-pull-back of the line bundle $\mathcal{I}$ and $s\equiv\mathbb{D}(|x|^\alpha \ff s_0)$ define a homogeneous $\mathbb{Z}/2$ harmonic spinor on $\mathbb{R}^3$.

	With the preceding understood, the rest of this paper constructs data $(Z, \mathcal{I}, \ff)$ on $S^2$ of the desired sort:  What is denoted by $Z$ is a set of $2N$ distinct points in $S^2$; what is denoted by $\mathcal{I}$ is a real line bundle defined in the complement of $Z$ with $Z$ being its points of discontinuity, and what is denoted by $\ff$ is a Laplace eigensection of $\mathcal{I}$ that obeys (1.3).  Note in this regard that condition on $|d\ff|$ in (1.3) is of paramount importance with regards to using this data set to construct a singularity model for $\mathbb{Z}/2$ harmonic 1-forms and spinors.  This condition on $d\ff$ is the only truly subtle issue.  

\section{Energy minimizing characterization}
	Let $\mathcal{T}$ denote the set of smooth sections of $\mathcal{I}$ over $S^2-Z$ subject to two constraints: 
\begin{itemize}
\item $\int_{S^2}|\ff^2|= 1$.
\item $|\ff|$ and $|d\ff|$ extend over $Z$ as H\"older continuous functions on $S^2$ that vanish on $Z$.
\end{itemize}
\hfill (2.1)

Supposing that $\ff \in \mathcal{T}$, define its `energy' to be the integral of $|d\ff|^2$:
$$\mathcal{E}(\ff) \equiv \int_{S^2} |d\ff|^2.$$
\hfill (2.2)\\
One might hope to find an eigensection of $\Delta$ that obeys (1.3) by minimizing the function $\mathcal{E}$ over the set $\mathcal{T}$.  (A minimizer of $\mathcal{E}$ is formally an eigensection.)  Unfortunately, there is no guarantee that $\mathcal{E}$ has a minimum in $\mathcal{T}$.   This is to say that minimizing sequences in $\mathcal{T}$ might converge to something that is not in $\mathcal{T}$.   The upcoming Proposition 2.1 makes a formal assertion to this effect.  To set the stage, introduce the set $\mathcal{T}_*$ consisting of smooth sections of $\mathcal{I}$ over $S^2-Z$ that obey the following:
\begin{itemize}
\item $\int_{S^2} |\ff|^2 = 1.$
\item $|\ff|$ extends over $Z$ as H\"older continuous function on $S^2$ that vanishes on $Z$.
\item $\int_{S^2} |d\ff|^2 < \infty.$
\end{itemize}
\hfill (2.3)\\
Notice the weaker condition on $|d\ff|$.  In particular, $\mathcal{T}\subset\mathcal{T}_*$.  Now, supposing that $\ff\in\mathcal{T}_*$, define $\mathcal{E}(\ff)$ as in (2.2).  

\begin{prop}
The infimum of $\mathcal{E}$ on $\mathcal{T}$ is the same as its infimum on $\mathcal{T}_*$.  Meanwhile $\mathcal{E}$ on $\mathcal{T}_*$ does take on its infimum value; and any section in $\mathcal{T}_*$ with this infimum value of $\mathcal{E}$ is an eigensection for the Laplacian.  Moreover any minimizing sequence in $T_*$ for $\mathcal{E}$ has a subsequence that converges to some minimizer of $\mathcal{E}$ in $\mathcal{T}_*$ as follows:  Let $\{f_n\}_{n\in \mathbb{N}}$ denote the subsequence (relabeled consecutively from 1) and let $\ff_*$ denote the corresponding minimizer of $\mathcal{E}$.  Then 
$$\lim_{n\to\infty}  \int_{S^2} (|d(\ff_n - \ff_*)|^2 + |\ff_n - \ff_*|^2) = 0.$$
\end{prop}
\noindent
This proposition is proved in Section 3.

	Our plan for circumventing this proposition is to use symmetry considerations.  To this end, let $G$ denote the group of orientation preserving symmetries of the regular tetrahedron.  This group can be viewed as a subgroup of $SO(3)$ by centering the tetrahedron so that the lines through the vertices and midpoints of the opposite edges intersect at the origin.  The vertices are taken to be the four points on $S^2$ with Euclidean coordinates 
	\vskip -1.5mm
$$(0,0,1), \left(-{\sqrt{2}\over 3},{\sqrt{2}\over \sqrt{3}}, -{1\over 3}\right), 
\left(-{\sqrt{2}\over 3},-{\sqrt{2}\over \sqrt{3}}, -{1\over 3}\right), 
\left({2\sqrt{2}\over 3}, 0, -{1\over 3}\right), $$
	\vskip -0.5mm
\hfill (2.4)\\
These points are labeled $p_1, p_2, p_3$ and $p_4$.  The group $G$ has a corresponding set of generators (subject to certain relations) denoted by $\{a_1, a_2, a_3, a_4\}$ with any given $a_k$ inducing a ${2\pi \over 3}$  rotation in the clockwise direction about the oriented line from the origin to the point $p_k$.  The composition $a_ia_j$  for $i\neq j$ is a $\pi$ rotation about the line through the origin and the midpoint of the edge of the tetrahedron that contains both $a_i$ and $a_j$ (thus, $a_ia_j = a_ja_i$).  Since this line also goes through the midpoint of the one edge that doesn't contain either $a_i$ or $a_j$, there are only three of these sorts of rotations in all.  (The group $G$ has 12 elements.)
Let $G_0$ denote the product group $\{1, -1\}  \times G$.  As indicated by the diagram in Figure 1 and explained subsequently, the group $G_0$ acts on the line bundle $\mathcal{I}$ as a group of isometries.
\begin{figure}[H]
 \centering
 \includegraphics[width=\textwidth]{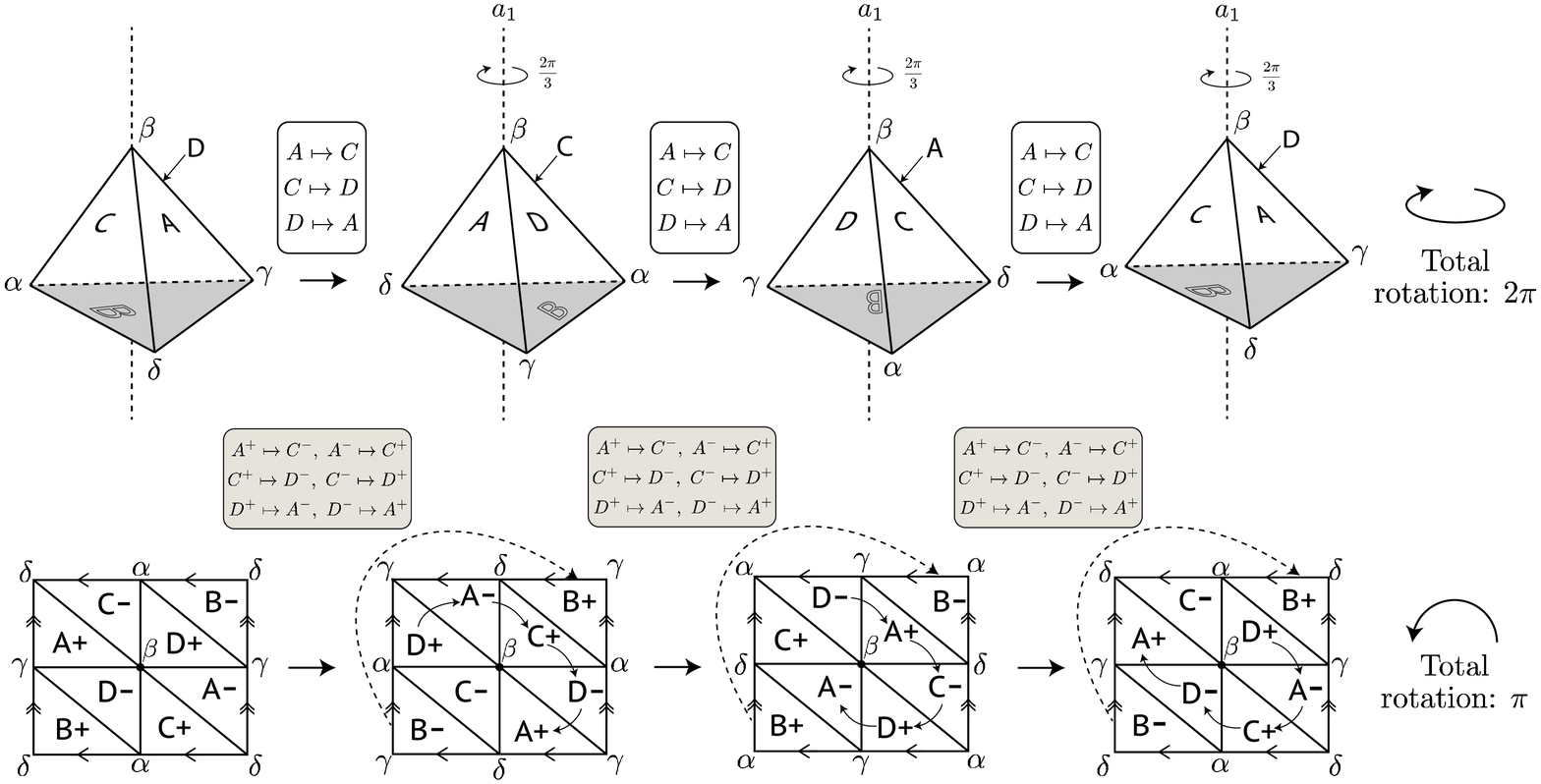}
 \caption{Rotations in $\{1, -1\} \times G$}\label{fig:rotation}
	\vskip -3mm
\end{figure}
To explain the diagram:  The tetrahedron maps homeomorphically to the sphere $S^2$ by the map $x \mapsto {x\over |x|}$.  This map is used to identify the sphere $S^2$ with the inscribed tetrahedron.  (Of particular note is that this identification respects the action of the tetrahedral subgroup of $SO(3)$.)  With the preceding identification understood, the set of points in $\mathcal{I}$ with norm 1 can be viewed as the complement of the branch points in the 2-fold branched cover of the tetrahedron branched over its vertices.  This 2-fold branched cover is a torus.  This torus is depicted in Figure 1 as a square with top and bottom edges identified and also with left and right-most edges identified.  With this depiction of the torus understood, then the left-most tetrahedron in Figure 1 and the left-most square in Figure 1 depict the branched covering map.  To make this map explicit, we have labeled the faces of the tetrahedron by   $A$,  $B$,  $C$  and  $D$  with opposite vertices labeled by $\alpha$, $\beta$, $\gamma$ and $\delta$.  The triangles in the torus are labeled by   $A$,  $B$,  $C$  and  $D$  (with an additional $\pm$ label).  The quotient branched covering map from the torus to $S^2$ sends the two triangles in the left-most torus in Figure 1 with any given letter to the like-labeled faces of the left-most tetrahedron in Figure 1.  (For example, the triangles labeled $A+$ and $A-$ in the torus are both mapped to the  $A$  labeled face in the tetrahedron.)  By the same token, the point in the torus labeled as $\beta$ (the middle point in the left-most torus of Figure 1) and the $\alpha$, $\gamma$, $\delta$ labeled points in the torus are mapped to the like labeled vertices in the tetrahedron by the branched covering quotient map.  These $\alpha$, $\beta$, $\gamma$, and $\delta$ labeled points in the torus are mapped to their counterparts in the tetrahedron in a 1-to-1 correspondence because they are the branching points of the covering.  

		The action of the product group $G_0 = \{1, -1\} \times G$ on the torus is depicted by Figure 1 using the other three vertical pairs of tetrahedron and torus.  To explain the depiction:  The element $a_1$ in the tetrahedral group $G$ effects a ${2\pi\over 3}$ clockwise rotation around the $\beta$ vertex moving face $D$ to $A$,  $A$  to $C$, and $C$ to $D$.  This rotation moves face $B$ to itself but rotates it by ${2\pi\over 3}$ in the process.  The group element $(-1,a_1)$ in the product group $G_0$ rotates each $D$ labeled triangle in the torus to its neighboring  $A$  labeled triangle; it moves each  $A$  labeled triangle to its neighboring $C$ labeled triangle, and it moves each $C$ labeled triangle to its neighboring $D$ labeled triangle.  It also switches the two $B$ labeled triangles (it must do this to preserve continuity).  The right-most three vertical pairs of tetrahedron and torus in Figure 1 illustrate the phenomena that 3 rotations of ${2\pi\over 3}$ around vertex $\beta$ of the tetrahedron results in a $180^\circ$ rotation around the point $\beta$ in the torus, a map that sends each letter labeled triangle with the additional label $+$ to the triangle with the same letter label and additional label $-$, and vice-versa.   The latter map of the torus (taking $\pm$ labels to $\mp$ labels with no change of letter label) gives the action of the element ($-1$, identity) in $G_0$.  Thus, the element ($-1$, identity) acts as the deck transformation element for the 2-fold covering map.  With regards to the action of the subgroup $\{1\} \times G$ of $G_0$:  This can be derived using Figure 1 by virtue of the fact that the elements $(1, a_1)$, $(1, a_2)$, $(1, a_3)$ and $(1, a_4)$ in this subgroup are the respective fourth powers of the elements $(-1,a_1)$, $(-1,a_2)$, $(-1,a_3)$ and $(-1,a_4)$. 
		
		Since the product group $G_0 = \{1, -1\}\times G$ acts isometrically on the bundle $\mathcal{I}$, it acts on the sets $\mathcal{T}$ and $\mathcal{T}_*$ because both sets of conditions are described by rotationally invariant conditions.  With the preceding understood, let $\mathcal{T}_0$ denote the subset of sections in $\mathcal{T}$ where the elements $\{(-1,a_1), (-1,a_2), (-1,a_3), (-1,a_4)\}$ each acts as multiplication by $-1$.  The $\pm$ pattern of labels in Figure 1 for the torus triangles indicates that this set $\mathcal{T}_0$ is not empty.  To see why, view a section of $\mathcal{I}$ as $\mathbb{R}$-valued function defined on the complement in the torus of the four points $\{\alpha, \beta, \gamma, \delta\}$ that is equivariant with respect to the action of the deck transformation element ($-1$, identity) in $G_0$.  (A function is equivariant with respect to this action if and only if it has opposite signs on any two triangles with the same letter label.)  With sections viewed in this light as functions, the pattern of $\pm$ signs in Figure 1 indicates the sign of an equivariant function of the desired sort that vanishes on the edges of the triangle.  Indeed, the action of any element from the four-element subset $\{(-1,a_1), (-1,a_2), (-1,a_3), (-1,a_4)\}$ changes the sign of a function with the indicated pattern of signs because each element from this set moves any given triangle to an adjacent one, and because no two adjacent triangles in the torus have the same sign. 
\begin{prop}
The function $\mathcal{E}$ on $\mathcal{T}_0$ takes on on its infimum value; any section in $\mathcal{T}_0$ with this infimum value of $\mathcal{E}$ on $\mathcal{T}_0$ is an eigensection for the Laplacian that obeys the conditions in (1.3).
\end{prop}

As we show in Section 4, an eigensection in $\mathcal{T}_0$ for $\Delta$ behaves near each point of $p\in Z$ with respect to a complex coordinate centered at that point (call it $u$) as the real part of a non-zero complex multiple of either $u^{k+{1\over2}}$ or $\bar u^{k+{1\over2}}$  plus terms that are $o(u^{k+{1\over2}})$ with $k$ being a positive integer.  (So $k \geq 1$.)  We also show that the differential of the section vanishes near $p$ as $|u|^{k-{1\over2}}$.

	As a parenthetical remark, the arguments for Proposition 2.2 can be used with only minor changes to prove that there is an infinite set of normalized eigensections of the Laplacian in the space $\mathcal{T}_0$ such that the values of the function $\mathcal{E}$ diverges along any sequence with no convergent subsequence.  Each of these supplies a homogeneous $\mathbb{Z}/2$ harmonic 1-form on $\mathbb{R}^3$ (but all have the same set $Z$, the four rays from the origin to the vertices of the inscribed tetrahedron.)

\section{Proof of Proposition 2.1}
	The proof has five parts.\\

	\it Part 1:  \rm This part of the proof explains why the infimums of $\mathcal{E}$ on $\mathcal{T}$ and $\mathcal{T}_*$ are the same.   To this end, fix once and for all a smooth, non-increasing function on $\mathbb{R}$ to be denoted by $\chi$ that equals 1 for $t<{1\over 4}$  and equals 0 for $t \geq {3\over 4}$.   With $\chi$ in hand, fix a positive number to be called $\rho$ with its upper bound being ${1\over 1000}$  times the minimum of the distances between the points in $Z$.   Set $\chi_\rho$ to be the function on $S^2$ that is given by the rule $\displaystyle\chi_\rho(\cdot)=\prod_{p\in Z}\chi\left(2-{\text{dist}(\cdot, p)\over \rho}\right)$. This function is equal to 1 where the distance to $p$ is greater than $2\rho$, and it is equal to zero where the distance to $Z$ is less than $\rho$.  Note that the derivative of $\chi_\rho$ is non-zero only in those annuli with inner radius $\rho$ and outer radius 2$\rho$ centered at the points in $Z$.   Note also:  The norm of the derivative of $\chi_\rho$ is at most some $\rho$-independent multiple of $1\over \rho$.   This implies that the integral of  $|d\chi_\rho|^2$ is bounded as $\rho\to0$ because the area where $d\chi_\rho$ is non-zero is bounded by a $\rho$-independent multiple of $\rho^2$.
	
	With $\chi_\rho$ understood, and supposing that $\ff$ is an element in $\mathcal{T}_*$, then $\chi_\rho\ff$ is an element in $\mathcal{T}$ because it is zero where the distance to $Z$ is less than $\rho$.  Because  $|\ff|$  is zero on $Z$ and is H\"older continuous, the integral of  $|\chi_\rho\ff|^2$ over $S^2$  differs from that of  $|\ff|^2$ by at most an $o(1)$ multiple of $\rho^2$ when $\rho$ is small.  This is because  $|\chi_\rho\ff|$  differs from  $|\ff|$  only where $\chi_\rho\neq 1$, which is a union of $N$ disks in $S^2$ each with area at most $4\pi\rho^2$ and because $\ff$ in each of these disks has very small norm.  
	
Meanwhile:  The integral of  $|d(\chi_\rho\ff)|^2$ over $S^2$  differs from that of  $|d\ff|^2$ by a small number that also limits to zero as $\rho\to0$.  Indeed, the difference is bounded by a $\rho$ independent multiple of the sum of two integrals.  The first is the integral of  $|d\ff|^2$ over the union of the $N$ disks of radius 2$\rho$ where $\chi_\rho\neq 1$; and the latter integral has limit zero as $\rho\to0$ because the area of the integration domain goes to zero in this limit.  The second integral is that of  $|d\chi_\rho||\ff|^2$ over these same disks.  That integral also has limit zero as $\rho\to0$.  Although the $\rho\to0$ limit of the integral of  $|d\chi_\rho|^2$ is not zero, the integral of  $|d\chi_\rho||\ff|^2$ limits to zero as $\rho$ does because the maximum of  $|\ff|$  where $d\chi_\rho\neq 0$ limits to zero as $\rho$ does (by virtue of the fact that  $|\ff|$  is zero on $Z$). 

	The remarks in the preceding two paragraphs imply directly that the infimums of $\mathcal{E}$ on $\mathcal{T}$and $\mathcal{T}_*$ are the same.  Indeed, if $\epsilon > 0$ and if $\ff\in\mathcal{T}_*$ is such that $\mathcal{E}(\ff)$ is less than $\epsilon$ above the infimum of $\mathcal{E}$ on $\mathcal{T}_*$, then the preceding paragraphs imply that there is a number $m_\rho$ defined for positive but very small $\rho$ such that $m_\rho\chi_\rho\ff\in\mathcal{T}$ and $\mathcal{E}(m_\rho\chi_\rho\ff) <2\epsilon$.\\
	
	\it Part 2:  \rm To find a minimizer for $\mathcal{E}$ in $\mathcal{T}_*$, it proves useful along the way to introduce a Hilbert space to be denoted by $L^2_1(\mathcal{I})$, which is defined as follows:  It is the completion of the vector space of real number multiples of elements in $\mathcal{T}_*$ using the norm whose square is given by the rule
	$$\frak{f} \to \int_{S^2}(|d\ff|^2 + |\ff|^2).$$
	\hfill (3.1)\\
	This Hilbert space norm of $\mathfrak{f}$ is denoted here by $\|\ff\|$.  Note that $\|\ff\|^2 = \mathcal{E}(\ff) + 1$ when the integral of $|\ff|^2$ on $S^2$ is equal to 1.  Thus, an a priori bound on $\mathcal{E}(\ff)$ gives an a priori bound on $\|\ff\|^2$ when the integral of $|\ff|^2$ is 1.  This is the motivation for introducing this norm and the associated Hilbert space $L_1^2(\mathcal{I})$.
	
	The lemmas that follow list some basic facts about $L_1^2(\mathcal{I})$.  The first lemma uses $\langle v, w\rangle$ to denote the Euclidean metric inner product between 1-forms $v$ and $w$.
\begin{lemma}
 If $\{\ff_n\}_{n\in\mathbb{N}}$ is a sequence in $L_1^2(\mathcal{I})$ with bounded norm, then there is an element $\ff_* \in L_1^2(\mathcal{I})$ and a subsequence $\Lambda \subset \mathbb{N}$ such that the following are true:
 \begin{itemize}
\item	If $\ff^\prime$ is any given element in $L_1^2(\mathcal{I})$, then $\displaystyle\lim_{n\in \Lambda} \int_{S^2} \langle d(\ff_n - \ff_*), d\ff^\prime\rangle = 0.$
\item $\displaystyle \int_{S^2}|d\ff_*|^2 \leq \liminf_{n\in \Lambda} \int_{S^2} |d\ff_n|^2$.
\item $\displaystyle \lim_{n\in\Lambda} \int_{S^2} | \ff_n - \ff_*|^2 = 0 $.
\end{itemize}
\end{lemma}

\begin{lemma}
Fix a positive number (to be called $\rho$) that is less than ${1\over 4}$  times the minimum distance between any two points in $Z$; then fix a non-negative $\delta < \rho$.  Let $A$ denote an annulus in $S^2$ centered about a point in $Z$ with inner radius $\delta$ and outer radius $\rho$.  If $\ff\in L_1^2(\mathcal{I})$, then $\displaystyle\int_A|\ff|^2 \leq 4\rho^2\int_A|d\ff|^2$.
\end{lemma}

The rest of this part of the proof is occupied with the proofs of these two lemmas.\\

\noindent \emph{Proof of Lemma 3.1}:  \rm Before starting, introduce the nested sequence $\{\Omega_m\}_{m\in \mathbb{N}}$ of subspaces in $S^2-Z$ with any given $\Omega_m$ denoting the set of points which have distance ${1\over m}$ or more from each point in $Z$.   The union of all of these is the whole of $S^2-Z$. 

	To prove the assertion of the top bullet of the lemma, it is sufficient to observe that the Hilbert space $L^2_1(\mathcal{I})$ has a countable dense set, which is to say that it is separable (the Banach-Alaoglu theorem).   This fact about a countable dense set follows by virtue of two facts:  First, the set of sections of $\mathcal{I}$ that vanish outside of any given $m \in \mathbb{N}$ version of $\Omega_m$ has a countable dense set with respect to the $\|\cdot\|$-norm topology.  (A countable basis consists of eigensections of $\Delta$ that are zero on the boundary of $\Omega_m$.)  Second, the collection $\{\Omega_m\}_{m\in\mathbb{N}}$ is observably countable.  The second bullet’s assertion also follows from the Banach-Alaoglu theorem.  
	
	To prove the assertion of the third bullet, note first that the sequence $\{|\ff_n|\}_{n\in \Lambda}$ is a bounded sequence in the $L^2_1$ Sobolev space of functions on $S^2$.  It therefore has a weakly convergent subsequence.  The latter converges strongly to its limit in the $L^2$ topology by virtue of the fact that the forgetful map from $L^2_1$ to $L^2$ is compact.  The limit in $L^2$ is necessarily $|\ff_*|$.   As a consequence, 
$$\int_{S^2} |\ff_*|^2 = \lim_{n\in \Lambda}\int_{S^2}|\ff_n|^2.$$
\hfill (3.2)\\
This last inequality implies that the $n \to \infty$ limit (for $n \in\Lambda$) of the term in parenthesis on the right-hand side of the next identity is zero.
$$\int_{S^2}|\ff_n - \ff_*|^2 = \left(\int_{S^2}|\ff_n|^2 - \int_{S^2}|\ff_*|^2\right) - 2\int_{S^2}(\ff_n - \ff_*)\ff_*.$$
\hfill (3.3)\\
That conclusion implies in turn that $n\to \infty$ limit (for $n\in\Lambda$) of the left-hand side of (3.3) is zero because the $n\to \infty$ limit (for $n\in\Lambda$) of the right-most term in (3.3) is zero due to the weak convergence of $\{\ff_n\}_{n\in\mathbb{N}}$ to $\ff_*$.\\

\emph{Proof of Lemma 3.2}:  The disk of radius $2\rho$ centered at a point in $Z$ has a radial coordinate $r$ and angle coordinate $\theta$ wherein the metric has the form $dr^2 + \sin^2 r d\theta^2$.  The integral of  $|d\ff|^2$ over the annulus $A$ with inner radius $\delta$ and outer radius $\rho$ centered at that point can be written using these coordinates as
$$\int_\delta^\rho\left(\int_0^{2\pi} \left(
\left|{\partial \over \partial r}\ff\right|^2
+
{1\over \sin^2 r}\left|{\partial \over \partial \theta}\ff\right|^2 
\right)d\theta\right)\sin r dr.$$
\hfill (3.4)\\
The restriction of $\mathcal{I}$ to any constant $r$ circle is the M\"obius line bundle (the unorientable line bundle over $S^1$).  Keeping in mind that the smallest eigenvalue of $-{\partial^2\over \partial \theta^2}$  acting on sections of that line bundle is ${1\over 4}$, it follows that the integral in (3.4) is no smaller than
$$\int_\delta^\rho\left(\int_0^{2\pi} 
{1\over 4\sin^2 r}\left|\ff\right|^2
d\theta\right)\sin r dr.$$
\hfill (3.5)\\
That in turn is no smaller than  ${1\over 4\rho^2}$ times the integral of  $|\ff|^2$ over $A$.  This last observation leads directly to what is asserted by the lemma. \\
 
\it Part 3:  \rm This part and Parts 4 and 5 of the proof explain why there is necessarily a section $\ff\in\mathcal{T}_*$ with $\mathcal{E}(\ff)$ being the infimum of $\mathcal{E}$ on $\mathcal{T}_*$.  To this end, let $\mathcal{E}_*$ denote the number $\inf_{f\in \mathcal{T}_*}\mathcal{E}(\ff)$, the infinium of $\mathcal{E}$ on $\mathcal{T}_*$.  Fix a minimizing sequence $\{\ff_n\}_{n\in\mathbb{N}} \in \mathcal{T}_*$ such that $\displaystyle \lim_{n\to\infty} \mathcal{E}(\ff_n)= \mathcal{E}_*$. By weeding out members if necessary (and subsequently renumbering consecutively from 1), we can assume that $2\mathcal{E}_* \geq \mathcal{E}(\ff_n) \geq \mathcal{E}(\ff_{n+1})$ for all $n \geq 1$.  Lemma 3.1 supplies a subsequence $\Lambda \subset\mathbb{N}$ and an element $\ff_*\subset L_1^2(\mathcal{I})$ that is described by the three bullets of that lemma.  Note in particular that the integral of  $|\ff_*|^2$ over $S^2$  is equal to 1 by virtue of the third bullet in Lemma 3.1.  This implies (among other things) that $\ff_*$ is not identically zero.   

As explained next, $\mathcal{E}(\ff_*)$ is equal $\mathcal{E}_*$. To explain why this is, note first that any element in $L_1^2(\mathcal{I})$ is (by definition) the norm convergent limit of some sequence from $\mathcal{T}_*$.  In particular, this is true of $\ff_*$:  There is some sequence, call it $\{\mathfrak{g}_m\}_{m\in\mathbb{N}}$, such that $\lim_{m\to\infty} \|\mathfrak{g}_m - \ff_* \|  = 0$.  By virtue of this (and because the function $\mathcal{E}$ appears as part of the norm), it follows that   
$$\int_{S^2} |d\ff_*|^2 = \lim_{m\to\infty} \mathcal{E}(\mathfrak{g}_m).$$
\hfill (3.6) \\
This implies in turn that $\mathcal{E}(\ff_*) \geq \mathcal{E}_*$ (because this is the case for all $\mathcal{E}(\mathfrak{g}_m))$.  Then, by virtue of the  second bullet of Lemma 3.1, it follows that $\mathcal{E}(\ff_*)$ must be equal to $\mathcal{E}_*$.   With it understood that $\mathcal{E}(\ff_*) = \mathcal{E}_*$, fix $n\in\Lambda$ for the moment and write that
$$\int_{S^2} |d( \ff_n -  \ff_*)|^2 = \int_{S^2}|d \ff_n|^2 - 2 \int_{S^2} \langle d( \ff_n -  \ff_*), d \ff_*\rangle - \int_{S^2}|d \ff_*|^2.$$
\hfill (3.7)\\
The $ \ff^\prime = \ff_*$ instance of the top bullet of Lemma 3.1 says that the middle term of the right-hand side in (3.7) has limit zero as $n\to\infty$ (for $n\in \Lambda$). Meanwhile, the right-most term of the right-hand side of (3.7) is $\mathcal{E}_*$; and the left-most term on the right-hand side of (3.7) limits to $\mathcal{E}_*$ as $n\to\infty$.   As a consequence, the $n\to\infty$ limit (for $n\in\Lambda)$ of the integrals on the left-hand side of (3.7) must vanish. \\

\it Part 4: \rm It remains at this point only to verify that $\ff_*$ is actually in $\mathcal{T}_*$.  The first point to make is that $\ff_*$ is smooth and that it is an eigensection for the Laplacian.  To this end, let $\mathfrak{u}$ denote an element from $\mathcal{T}_*$ with support in a small radius disk that is disjoint from the points in $Z$.  If the norm of $t\in\mathbb{R}$ is small enough (but still not zero), then $\ff_* + t\mathfrak{u}$ will be somewhere non-zero.  For such $t$, define $\mathcal{Z}(t)$ by the rule  
$$\mathcal{Z}(t)  =   \int_{S^2} |\ff + t \mathfrak{u}|^2.$$
\hfill (3.8)\\
Then, by virtue of the definition of $\mathcal{Z}$, the section $\mathcal{Z}(t)^{-{1\over2}}(\ff_*+t\mathfrak{u})$ is in $\mathcal{T}_*$.  And, as a consequence, the value of $\mathcal{E}(\cdot)$ on $\mathcal{Z}(t)^{-{1\over2}}(\ff_*+t\mathfrak{u})$ can never be less than $\mathcal{E}_*$.  This implies that the function $t \to \mathcal{E}(\mathcal{Z}(t)^{-{1\over2}}(\ff_*+t\mathfrak{u}))$ (which is defined for $t$ near 0) has a local minimum at $t=0$, because $\mathcal{E}_* = \mathcal{E}(\ff_*)$. And, that can happen only if   
$$\int_{S^2}\langle d\ff_*, d\mathfrak{u}\rangle - \mathcal{E}_*\int_{S^2}\ff_* \mathfrak{u} = 0$$
\hfill (3.9)\\
because the left-hand side of (3.9) is the first order Taylor approximation to the function of $t \to \mathcal{E}(\mathcal{Z}(t)^{-{1\over2}}(\ff_*+t\mathfrak{u}))$ at $t= 0$.  Since the condition in (3.9) holds for all section $\mathfrak{u}$, it follows using standard properties of the Laplacian on small radius disks disjoint from $Z$ (where $\mathcal{I}$ is isomorphic to the product $\mathbb{R}$ bundle) that $\ff_*$ is smooth and an eigensection for the Laplacian with eigenvalue $\mathcal{E}_*$. \\

\it Part 5: \rm This last part of the proof explains why  $|\ff_*|$  must vanish at the points in $Z$ and why  $|\ff_*|$  is uniformly H\"older continuous on a neighborhood of any such point.  There are eight steps to the explanation.  (The explanation follows arguments from Chapter 3.5 of \cite{morrey2009multiple}.)\\  

\underline{Step 1}:  Fix a disk centered at a point $p$ from $Z$ whose radius is much less than 1 and much less than the distance from $p$ to any other point in $Z$.   Supposing that $\rho$ is positive but less than half the radius of this disk, use the function $\chi$ from Part 1 to define a function on $S^2$ to be denoted by $\beta_\rho$ by the rule $\beta_\rho(\cdot)\equiv \chi\left(2{\text{dist}(\cdot, p)\over \rho}-1\right).$ This function is equal to 1 where the distance to $p$ is less than  ${1\over 2}\rho$ and it is equal to 0 where the distance to $p$ is greater than $\rho$.  Note that its derivative has support only in the annulus where the distance to $p$ is between  ${1\over 2}\rho$ and $\rho$; and that the norm of this derivative is bounded by a $\rho$-independent constant times  ${1\over \rho}$.  Next, given a positive number $\epsilon<\rho$, define a second function using $\chi$ to be denoted by $\mu_\epsilon$ by the rule $\mu_\epsilon(\cdot) \equiv \chi\left(2 \left(1-{\text{dist}(\cdot, p)\over \epsilon}\right)\right)$. This function is equal to 1 where the distance to $p$ is greater than $\epsilon$ and it is equal to 0 where the distance to $p$ is less than  ${1\over 2}\epsilon$.  Let $\mathfrak{u}_{\rho,\epsilon} =  \beta_\rho^2 \mu_\epsilon^2 \ff_*$.\\

\underline{Step 2}:  Let $D_{\rho \over 2}$ denote the disk of radius  ${\rho\over 2}$ centered at $p$; let $D_\epsilon$ denote the disk centered at $p$ with radius $\epsilon$; and let $A_\rho$ denote the annulus centered at $p$ with inner radius  ${1\over 2}\rho$ and outer radius $\rho$.  The following inequality holds by virtue of the $\mathfrak{u} = \mathfrak{u}_{\rho,\epsilon}$ version of (3.9), and by virtue of Lemma 3.2:
$$\int_{D_{\rho \over 2}} |d\ff_*|^2 \leq c {1\over \rho^2}\int_{A_\rho} |\ff_*|^2 + c{1\over \epsilon^2}\int_{D_\epsilon} |\ff_*|^2 + c\mathcal{E}_* \int_{D_{\rho \over 2}} |\ff_*|^2 
$$
with $c$ depending only on the choice of $\chi$.  In particular, it is independent of $\rho, \epsilon$, and $\ff_*$.  (This inequality is obtained by applying versions of the triangle inequality to terms that have derivatives of either $\beta_\rho$ or $\mu_{\epsilon}$. )\\
right most
\underline{Step 3}:  The preceding inequality is exploited by first invoking Lemma 3.2 to bound the three integrals on its right-hand side with the result being: 
$$(1 - 4c\mathcal{E}_*\rho^2) \int_{D_{\rho \over 2}}
|d\ff_*|^2 
\leq 
4c \int_{A_\rho} |d\ff_*|^2 + 4c\int_{D_\epsilon} |d\ff_*|^2.$$
 \hfill (3.11) \\
 To exploit this, suppose henceforth that $\rho$ is no greater than  ${1\over 2}(4c_0\mathcal{E}_*)^{-{1\over2}}$.   Assuming this, and seeing as the $\epsilon \to 0$ limit of the $D_\epsilon$ integral of  $|d\ff_*|^2$ is zero, taking $\epsilon$ ever small with limit zero leads from (3.11) to this:
$$\int_{D_{\rho \over 2}}
|d\ff_*|^2 
\leq 
8c \int_{A_\rho} |d\ff_*|^2.$$
\hfill (3.12)\\
 And, since the $A_\rho$ integral of  $|d\ff_*|^2$ is the difference between its $D_\rho$ and $D_{\rho \over 2}$ integrals, what is written in (3.12) leads (after rearranging) to: 
 $$(1+8c)\int_{D_{\rho \over 2}}
|d\ff_*|^2 
\leq 
8c \int_{D_\rho} |d\ff_*|^2.$$
\hfill (3.13)\\
This is to say that $$ \int_{D_{\rho \over 2}} |d\ff_*|^2\leq 
\gamma_0\int_{D_\rho} |d\ff_*|^2,$$
\hfill (3.14)\\
with $\gamma_0$ short hand for ${8c\over 1+8c}$ .  Of particular note is that $\gamma_0 < 1$.\\

\underline{Step 4}:  Now fix $\rho_0 > 0$ so that (3.14) holds with $\rho$ = $\rho_0$.  Then, for any positive integer $n$, let $\rho_n = 2^{-n}\rho_0$. Iteration of (3.14) (starting with $\rho_n$, then $\rho_{n-1}$, and so on to $\rho_1$) leads to the following: 
 $$ \int_{D_{2^{-n}\rho_0}} |d\ff_*|^2\leq 
\gamma_0^n\int_{D_{\rho_0}} |d\ff_*|^2.$$
\hfill(3.15) \\
Because $\gamma_0 < 1$, the preceding leads in turn to a bound of this sort:  If $\rho\in(0, \rho_*)$, then 
$$\int_{D_\rho}|d\ff_*|^2 \leq c_* \rho^{\alpha}$$
\hfill (3.16)\\ 
with $c_*$ being independent of $\rho$ and with $\alpha$ being the norm of the base 2 logarithm of $\gamma_0$:
$$\alpha = \left|\ln_2\gamma_0\right|.$$
\hfill (3.17) \\
(To obtain this, use (3.15) for that value of $n$ with the property that $2^{-n-1} \rho_0 < \rho < 2^{-n}\rho_0$.  Also, use the fact that the $D_{\rho_0}$  integral of  $|d\ff|^2$ is, in any event, no greater than $\mathcal{E}_*$.) 	 \\

\underline{Step 5}:  The inequality in (3.16) will now be used to define a value for  $|\ff_*|$  at the point $p$.  (Functions in the Sobolev space $L_1^2$ such as  $|\ff_*|$  can be ambiguous on a set of zero measure.)  To do this, fix a positive number $\rho < {1\over 4} \rho_0$ and having done this, define $f_\rho(p)$ to be the average of the function  $|\ff_*|$  on the radius $\rho$ circle centered at $p$.  Note in this regard that  $|\ff_*|$  is integrable on this circle by virtue of a standard Sobolev inequality (see Theorem 3.4.5 in \cite{morrey2009multiple}).  Now it follows via the fundamental theorem of calculus that if $\rho^\prime$ is positive, less than $\rho$ and greater than  ${1\over 4}\rho$, then 
$$|f_{\rho^\prime}(p) - f_\rho(p)| \leq c_1 \int_{D_\rho - D_{\rho^\prime}} {1\over \text{dist}(\cdot, p)}|d\ff_*|$$
\hfill (3.18)\\
with $c_1$ being independent of $\rho$ and $\rho^\prime$.  This inequality and (3.16) lead directly to the inequality 
$$|f_{\rho^\prime}(p) - f_\rho(p)| \leq c_2{\rho^\prime\over \rho}\rho^{\alpha\over 2}$$
\hfill (3.19)\\
with $c_2$ being independent of $\rho$ and $\rho^\prime$ (as long as $\rho^\prime$ is between  ${1\over 4}\rho$ and $\rho$).  

The inequality in (3.19) then implies that the function $\rho\to f_\rho$ converges uniformly as $\rho\to 0$  and that the approach to the limit (denote this limit by  $f_0$  for the moment) is such that   
$$|f_0 - f_\rho| \leq c_3 \rho^{\alpha \over 2}.$$
\hfill (3.20)\\
 To derive (3.20), iterate (3.21) by taking $\rho^\prime = {1\over 2}\rho$, then replacing $\rho$ by $\rho^\prime = {1\over 2}\rho$ and repeating and repeating and so on.  Then sum the resulting inequalities.  The sum converges because $\sum_{n\in \mathbb{N}}{1\over 2^{n\alpha}}$ converges. 
 
 With the preceding understood, the value of  $|\ff_*|$  at the point $p$ is \it defined \rm to be this number  $f_0$  (the $\rho\to 0$  limit of the $f_\rho$’s).  	\\

\underline{Step 6}:  Step away from $p$ for the moment to consider  $|\ff_*|$  at points near $p$ but not equal to $p$.  The purpose is to bound the variation of  $|\ff_*|$  in disks about a point $q \neq p$ with radius on the order of dist$(p,q)$ but less than dist$(p,q)$.  (It is assumed implicitly that dist$(p, q) <  {1\over 10}\rho_0$.)  The upcoming Step 7 proves the following:  Assume that $\rho_0 < 1$, that $\rho_0^2 \mathcal{E}_*$ is less than 1 and that nothing from $Z$ other than $p$ lie in the radius 2$\rho_0$ disk centered at $p$.  If $q$ has distance less than  ${1\over 10}\rho_0$ from $p$ and if $z$ has distance at most  ${1\over 8}$dist$(p,q)$ from the point $q$, then 
$$\bigg| |\ff_*(z)|  -  |\ff_*(q)|\bigg|\leq c_{*1} \text{dist}(p,q)^{\alpha\over 2}$$
\hfill (3.21)\\
with $c_{*1}$ being independent of $z$ and $q$. 

Granted (3.21), then what follows is a direct consequence:  If $\rho < {1\over 10}\rho_0$ and if $q$ is on the circle of radius $\rho$ centered at $p$, then  $|\ff_*|(q)$ differs from $f_\rho$ by at most $c_4 \rho^{\alpha\over 2}$ with $c_4$ being independent of $q$ and $\rho$.  This implies in turn (via (3.20)) that   
$$\big| |\ff_*| (p) -  |\ff_*|(q) \big| \leq c_5\text{dist}(q,p)^{\alpha \over 2}$$
\hfill (3.22)\\
with $c_5$ also independent of $p$ and $q$.  Thus, the function  $|\ff|(\cdot)$ is H\"older continuous at $p$.   

The fact that  $|\ff_*|$  is H\"older continuous at $p$ requires in turn that  $|\ff_*|$  must vanish at $p$.  The reason is as follows:  The line bundle $\mathcal{I}$ is non-trivial on all sufficiently small radius circles centered at $p$, and so $\ff_*$, being a section of $\mathcal{I}$, must vanish at one or more points on each of these circles.  In particular, there is a sequence of points converging to $p$ where  $|\ff_*|$  is zero.  Since  $|\ff_*|$  is H\"older continuous at $p$, it is continuous at $p$ and so  $|\ff_*|$  = 0 at $p$.\\ 

\underline{Step 7}:  This step explains why (3.21) holds.  To simplify notation in what follows, introduce $\upsigma$ to denote the distance between $p$ and $q$.  Let $D$ denote the disk in $S^2$  of radius  ${1\over 4}\sigma$ centered at $q$.  Because this disk is disjoint from p and from the rest of $Z$, the eigensection $\ff_*$ can be viewed as a real number valued function on $D$.  One can then use the Green’s function for the Laplacian on $D$ (with Dirichlet boundary conditions) to represent $\ff$ in $D$ as follows:  Let $z$ denote a point with distance less than ${1\over 8}\upsigma$ from $q$; and let $G_z$ denote the Green’s function for the Laplacian on the disk centered at $z$ with radius  ${1\over 8}\upsigma$ that vanishes on the boundary of this disk.  With regards to $G_z$:  There exists a number $c_6$ which is independent of $q$ and $z$ with the following significance:

\begin{enumerate}
\item	$\displaystyle|G_z| \leq c_6 \left|\ln \left({1\over \upsigma} \text{dist}(z,\cdot)\right)\right|$,
\item $\displaystyle|dG_z | \leq c_6{1\over \text{dist}(z,\cdot)}$,
\item $\displaystyle|\nabla dG_z|  \leq c_6{1\over \text{dist}(z,\cdot)^2}$.
\end{enumerate}
\hfill (3.23)\\
Let $\bar\omega_z$ denote the function $\chi\left({16\text{dist}(z,\cdot) \over \upsigma}-1\right)$.  Take $\mathfrak{u} = \bar \omega_z G_z$ in (3.9).  Since the integrands in (3.9) have support in the radius  ${1\over 8}\upsigma$ disk centered at $z$ (because of $\bar \omega_z$), an instance of integration by parts in the left-most integrals lead to the following identity:
$$ \ff_*(z) = \mathcal{E}_*\int_D \bar\omega_z G_z \ff_*
+ \int_D(-\Delta \bar\omega_z G_z + 2\langle d\bar \omega_z, dG_z\rangle)\ff_*.$$
\hfill (3.24)\\ 
What with the top bullet in (3.23), the absolute value of the left-most term on the right-hand side of (3.24) is necessarily bounded by  
$$ c_7\mathcal{E}_*\upsigma\left(\int_D |\ff_*|^2\right)^{1\over2}.$$
Here, $c_7$ is also independent of $z$ and $q$.  Meanwhile, the versions of the right-most term in (3.24) in the respective cases when $z = q$ and when $z$ is any point obeying dist$(z,q)< {1\over 8}\upsigma$ differ by at most  
$$c_8{1\over \upsigma}\left(\int_D |\ff_*|^2\right)^{1\over2}$$
\hfill (3.26)\\
with $c_8$ being independent of $z$ and $q$.   To prove this, use the lower two bullets in (3.23) to bound the derivative with respect to $z$ of the right-most term on the right-hand side of (3.24) by a factor of ${1 \over \upsigma}$  times what is written in (3.26). 

As for the integral of  $|\ff_*|^2$ in (3.25) and (3.26):  Its integral over $D$ is no larger than its integral over the disk of radius $2\upsigma$ centered at $p$ (not $q$) because $D$ is inside that larger disk.  The latter integral is no greater than $c_9\upsigma^{2+\alpha}$ by virtue of Lemma 3.2 and (3.16).  Here $c_9$ is independent of $\upsigma$.  Therefore, assuming that $\rho_0^2\mathcal{E}_*$ is less than 1, then the bounds in (3.25) and (3.26) lead directly to what is asserted in (3.21).  \\

\underline{Step 8}:  The previous steps proved that  $|\ff_*|$  is H\"older continuous at $p$.  This function is also H\"older continuous at points not in $Z$ because it is smooth away from $Z$, but this does not directly imply that  $|\ff_*|$  is \it uniformly \rm H\"older continuous near points in $Z$.  To see that it is, fix $p\in Z$ and suppose that $q$ is again a point in $S^2$  but not from $Z$.   Assume that $q$ has distance less than  ${1\over 100}\rho_0$ from $p$.  Let $q^\prime$ denote a second point from $S^2-Z$ which also has distance less than  ${1\over 100}\rho_0$ from $p$.  There are two cases to consider:  The first occurs when dist$(q^\prime,q) \geq dist(q,p)$ and the second when this inequality is violated.  In the first case, the inequality in (3.21) with $z = q^\prime$ implies that  
$$\big| | \ff_*(q^\prime)|  - | \ff_*(q)|\big| \leq c_{*1} 100^{\alpha\over2}\text{dist}(q^\prime,q)^{\alpha\over 2}.$$
\hfill (3.27)\\
In the second case, the distance between $q^\prime$ and $q$ is at most ${1\over 10}$  times the distance from either to $p$.  In this case, the identity in (3.24) can used with $z$ situated on the short geodesic arc between $q^\prime$ and $q$ (this arc is denoted by $I$).  Differentiating that identity with respect to $z$ and using (3.23) leads to the following bound for  $|d\ff|$  along the arc $I$:   
$$|d\ff_*|\leq c_9\left(\mathcal{E}_*\text{dist}(\cdot,p) + {1\over \text{dist}(\cdot, p)}\right)\sup_D|\ff_*|$$
\hfill (3.28) \\
with $c_9$ denoting a number that is independent of the point in question.   (The function $\dist(\cdot, p)$ appears here because the radius of the disk $D$ is on the order of $\dist(\cdot,p)$ when the point in question is arc $I$.) 

With regards to $\sup_D|\ff_*|$: It is bounded by $q$ and $q^\prime$ independent multiple of $\sigma^{\alpha\over2}$, that this is so follows from the versions of (3.22) with $q$ allowed to be any point in $D$. With this bound understood, then (3.28) has the following implications:  In the case when $\alpha\geq 2$, it implies that    
$$|\ff_*(q^\prime) - \ff_*(q)|  \leq c_{10} \left(\dist(q,p)^{{\alpha\over2}-1}+ \dist(q^\prime,p)^{{\alpha\over 2}-1}\right)\dist(q^\prime,q)$$
\hfill (3.29)\\
with $c_{10}$ being independent of $q$ and $q^\prime$.  (The distance from $z$ to $p$ for $z\in I$ is not less than the sum of the distances from $q$ to $p$ and from $q^\prime$ to $p$.)  In the case when $\alpha < 1$, one has   
$$|\ff_*(q^\prime) - \ff_*(q)|  \leq c_{11}\left(\sup_{z\in I}{1\over \dist(z,p)^{1-{\alpha\over 2}}} \right)\dist(q^\prime, q).$$
\hfill (3.30)\\
(The number $c_{11}$ is independent of $q$ and $q^\prime$.)  Write this last inequality as 
$$|\ff_*(q^\prime) - \ff_*(q)|  \leq c_{10} \left( \sup_{q\in I}  {\dist(q^\prime, q)^{1-{\alpha \over 2}}\over \dist(z, q)^{1-{\alpha \over 2}}}\right)\dist(q^\prime, q)^{\alpha \over 2}$$
\hfill (3.31)\\
to see that it implies in turn that   
$$|\ff_*(q^\prime) - \ff_*(q)| \leq c_{12} \dist(q^\prime,q)^{\alpha \over 2}$$
\hfill (3.32) \\
with $c_{12}$ being independent of $q$ and $q^\prime$.  This is because distance from $z$ to $p$ when $z\in I$ is no \it smaller \rm than the sum of the respective distances from $q$ to $p$ and from $q^\prime$ to $p$. 

This last inequality with (3.29) and (3.27) and (3.22) prove that  $|\ff|$  is uniformly H\"older continuous near $p$.  

\section{ Proof of Proposition 2.2} 	The proof of this proposition has six parts.  \\

\it	Part 1: \rm The binary tetrahedral group $G_0$ acts on $L^2_1(\mathcal{I})$ by isometries.  The theory of finite group representations (see, Chapter 3 in Mackey’s book \cite{mackey1978unitary}) leads to the following two observations:  First, $L^2_1(\mathcal{I})$ has an orthogonal (with respect to the Hilbert space norm), direct sum decomposition as $\lm^{-} \oplus \lm^{-\perp}$ with $\lm^{-}$ denoting the subspace of sections where the generators $\{\alpha_1, \alpha_2, \alpha_3, \alpha_4\}$ act as multiplication by $-1$.  Second, this decomposition is orthogonal with respect to the $L^2$ inner product also (the $L^2$ inner product comes from the norm whose square sends a section $\ff$ to the integral of  $|\ff|^2$ over $S^2$). \\

Granted the preceding two facts, then the arguments from Section 3 can be repeated almost verbatim but for $\lm^{-}$ replacing in each instance of $L^2_1(\mathcal{I})$ to see that  \\

\begin{itemize}
\item There exists a section $\ff_0$ in $\mathcal{T}_* \cap \lm^{-}$ that minimizes the function $\mathcal{E}$ on $\mathcal{T}_* \cap \lm^{-}$. 
\item Let $\mathcal{E}_0$ denote this minimal value of $\mathcal{E}$ on $\mathcal{T}_* \cap \lm^{-}$.  The section $\ff_0$ is a eigensection for the Laplacian on $\mathcal{T}_*$ with eigenvalue $\mathcal{E}_0$. 
\end{itemize}
\hfill (4.1)\\
 With regards to proving that $\ff_0$ is an eigensection:  The verbatim repeat of the arguments up to (3.9) find that (3.9) holds with $\ff_*$ replaced by the minimizer $\ff_0$ and $\mathcal{E}_*$ replaced by $\mathcal{E}_0$ and with $\mathfrak{u}$ 
\it restricted \rm to $\lm^{-}$.  However, (3.9) also holds when $\mathfrak{u}$ is orthogonal to $\lm^{-}$ because the decomposition $L^2_1(\mathcal{I}) = \lm^{-}\oplus\lm^{-\perp}$ is orthogonal for both the Hilbert space inner product and the $L^2$ inner product.  This is to say that both integrals that appear in (3.9) are zero as long as $\ff_*$ is from $\lm^{-}$ and $\mathfrak{u}$ is from $\lm^{-\perp}$ (no assumption is necessary in this regard about $\ff_*$ minimizing anything).   	\\
 
 \it Part 2: \rm  It remains now to prove that  $|d\ff_0|$  also extends over $Z$ as a H\"older continuous function (which is to say that $\ff_0$ is in the space $\mathcal{T}_0$).  To this end, fix a point in $Z$ to be denoted by $p$ and then introduce a stereographic coordinate centered at $p$ to identify the radius $\rho_0$ disk centered at $p$ with a small radius disk in $\r2$ centered at the origin.  Use $z$ to denote the complex coordinate on $\r2$; thus $z=0$ is the point $p$ and the disk of radius $\rho_0$ around $p$ is the $|z|  < \rho_1$ disk in $\mathbb{C}$ with $\rho_1 = \rho_0 + \mathcal{O}(\rho^3_0)$.  The function  $|z|$  is denoted by $r$ and the argument of $z$ is denoted by $\theta$ (with $\theta \in\real/2\pi\integer$).
 
    	The restriction of the function $\ff_0$ to any constant $r$ circle (with $r < \rho_1$) is a section of the M\"obius line bundle over the circle.  As such, it can be written as a linear combination of eigensections on this circle for the circle Laplacian.  Since the Laplacian on the circle is ${d^2\over d\theta^2}$, the corresponding set of Laplace eigensections for the M\"obius bundle is the collection $\left\{ e^{i(n + {1\over 2})\theta}\right\}_{n\in \integer}$. Thus, $\ff_0$ has the Fourier decomposition:  
	$$\ff_0|_{(r,\theta)} = \sum_{n\in \integer} \mathfrak{a}_n(r)e^{i\left(n + {1\over 2}\right)\theta},$$
\hfill (4.2)\\
with $\{\fa_n(\cdot)\}_{n\in\integer}$ denoting functions on $(0,\rho_1]$.

   	Let $\alpha$ denote the generator from the set $\{\alpha_1, \alpha_2,\alpha_3, \alpha_4\}$ that fixes the point $p$.  Because this generator acts on $S^2$  as the ${2\pi\over 3}$   rotation about the point $p$, it appears with respect to the coordinate $z$ as the ${2\pi\over 3}$  rotation about the $z=0$ point in $\mathbb{C}$ which is to say that its action fixes $r$ and sends $\theta$ to $\theta+{2\pi \over 3}$.  This action lifts to an action on the M\"obius line bundle on the fixed radius circles which sends any given $n\in\integer$ version of $e^{i\left(n+{1\over 2}\right)\theta}$  to the section $e^{{2n + 1 \over 3}i}e^{i\left(n+{1\over 2}\right)\theta}$. Thus, each Laplace eigenfunction is sent to a multiple of itself by this action.  But note that this multiple is equal to $-1$ if and only if $2n+1$ is an odd multiple of 3 which is to say that $n$ is congruent to 1 (mod 3).  Therefore, the expansion in (4.2) can be written as  
	$$\ff_0|_{(r,\theta)} = \sum_{m\in\integer} \fa_{3m+1}(r)e^{i3\left(m+{1\over 2}\right)\theta}$$
\hfill (4.3)\\

	Looking ahead, the absence in (4.3) of the $n = 0$ and $n = -1$ eigensections is the key input to the proof of Proposition 2.2.  \\
	
	\it Part 3: \rm By virtue of the fact that $\ff_0$ is an eigensection for the Laplacian on $S^2$, any given $\fa_n(\cdot)$ that appears in (4.3) (or in (4.2) if $\ff_0$ is not constrained with respect to $G_0$) must obey the differential equation  
	$$-r{d\over dr}r {d\over dr} \fa_n + \left(n+{1\over 2}\right)^2 \fa_n =  \mathcal{E}_* {4r^2\over (1+r^2)^2} \fa_n. $$
	\hfill (4.4)\\
	This equation is of \it Sturm-Liouville \rm type on the domain $[0, \rho_1]$ with 0 being a \it regular \rm point and all other points being \it ordinary. \rm  (See Chapters 10.2 and 10.3 of \cite{whittaker1996course} for the definitions of the italicized terms.)  What this implies in the case at hand (see Chapter 10.3 of \cite{whittaker1996course}) is that $\mathfrak{a}_n$ on $[0, \rho_1]$ can be written as   
	$$\fa_n =  r^{\left|n + {1\over 2}\right|}(1+\fu_n(r))$$
	\hfill  (4.5)\\
	 with $\fu_n$ being a real analytic function near $r = 0$.  (Only the positive exponent appears in (4.5) because  $|\ff_0|$  vanishes at $z=0$.) 
	 
	     	Keeping in mind that those $n$ that appear in (4.3) are congruent to 1 (mod 3), the prefactor powers of $r$ that appear in (4.5) are no less than ${3\over 2}$.  This suggests (strongly) that   
		$$|\ff_0| \leq \mathcal{O}\left(r^{3\over 2}\right) \quad and \ that \quad |d\ff_0| \leq \mathcal{O}\left(r^{1\over 2}\right).$$
		\hfill (4.6)\\
		 near $p$.  The subsequent parts of Proposition 2.2's proof explain why (4.6) is an accurate depiction of the behavior of  $|\ff|$  and  $|d\ff|$.  	\\
		 
		 \it Part 4: \rm To prove that (4.6) is accurate, fix for the moment a large, positive integer $N$ and use (4.3) to write $\ff_0$ where  $r\leq\rho_1$ as  
		 $$\ff_0 = \sum_{m\in\integer; \ \left|m+{1\over 2}\right|< N}\fa_{3m+1}(r)e^{i3\left(m + {1\over 2}\right)\theta} + \ff_N$$
		 \hfill (4.7)\\
 where $\ff_N$ is the sum of the terms in (4.3) with  $|m|  \geq N$.   The left-most term on the right-hand side of (4.7) is described by (4.6) because it is a finite sum of terms with each term having norm bounded by a contant multiple of $r^{3\over2}$ and with the norm of each term’s differential bounded by a constant multiple of $r^{1\over2}$.    
 
 	As for the $\ff_N$ part of (4.7):  The first point to note is that $\ff_N$ is an eigensection with eigenvalue $\mathcal{E}_*$ for the Laplace operator acting on sections of $\mathcal{I}$ over the disk where the coordinate $r$ is at most $\rho_1$. It is also the case that the  $L^2_1$ norm of $\ff_N$ on this disk is no greater than that of $\ff_0$ and likewise for its  $L^2$ norm.  This is because $\ff_N$ is orthogonal to the left-most sum on the right-hand side of (4.7) in both norms.   The key point with regards to $\ff_N$ is that it obeys a version of Lemma 3.2 with the number 4 replaced by a much smaller number: 
	
	\begin{lemma}
	Fix a positive number (to be called $\rho$) that is less than $\rho_0$; then fix a non-negative $\delta<\rho$.  Let $A$ denote an annulus in $S^2$ centered about the point $p$ with inner radius $\delta$ and outer radius $\rho$.  Then  
	$\displaystyle\int_A|\ff_N|^2 \leq {4\over (2N+1)^2} \rho^2 \int_A |d\ff_N|^2.$
	\end{lemma}
	
	\emph{Proof of Lemma 4.1}: \rm 
	Repeat the argument for Lemma 3.2 noting that the factor of ${1\over 4}$  that appears in (3.5) can be replaced by ${(2N+1)^2\over 4}$  since this is the norm of the smallest eigenvalue of $-{d^2\over d\theta^2}$  acting on  the relevant vector space of sections of the M\"obius bundle over the circle. \\
	
	\it  	Part 5: \rm Granted this lemma, then a repetition of the arguments in Steps 3 and 4 of in Part 5 of the previous section lead to the following analogy of (3.16):  
	$$\int_{D_\rho} |d\ff_N|^2 \leq c_*\rho^\alpha$$
	\hfill (4.8)\\
	 with $\alpha$ now given by   $$\alpha =  \left|\ln_2 \left({8c \over (2N+1)^2 + 8c}\right)\right|  $$
	 \hfill (4.9)\\
	  \Big(The effect of Lemma 4.1 is to replace the number $c$ by ${c \over (2N+1)^2 }$.\Big)  
	  
	  	One can now repeat Steps 5 and 6 of the Part 5 in the previous subsection to see that (3.22) holds with $\ff_N$ appearing instead of $\ff_*$ when $\dist(q,p)\leq {1\over 10}\rho_0$.  This version says   
	  $$|\ff_N |(q) \leq c_{13}\dist(q,p)^{\alpha \over 2}$$
\hfill  (4.10) 
\\
if $\dist(q,p) \leq{1\over 10} \rho_0$.  Likewise, the arguments from Steps 7 and 8 of the preceding subsection that lead to (3.28) can be repeated to rederive (3.28) which says (given (4.10)) that 
 $$|d\ff_N| (q) \leq c_{14} \dist(q,p)^{\alpha \over 2}$$
 \hfill   (4.11)\\
  if $\dist(p,q)  \leq {1\over 100}\rho_0$.  
  
   The preceding bound and (4.10) imply that (4.6) holds for $\ff_N$ and thus for $\ff_0$ if $N$ is sufficiently large.  	\\
  
  \it Part 6:  \rm  The proof that the function  $|d\ff_0|$  is uniformly H\"older continuous near any given point in $Z$ is much like the proof in Steps 7 and 8 of Part 5 of the previous section for the analogous assertion about  $|\ff_0|$.  The starting point for this is the identity in (3.24), which one differentiates to obtain an identity for $d\ff_*$.
  
     The details of the argument are left to the reader except for one comment which concerns the left-most term on the right-hand side of (3.24).  Comparing the derivative of this term at a point $q$ and at a nearby point $q^\prime$ is slightly subtle by virtue of the fact that the norm of the second derivative of $G_z$ is singular at $z$.  Circumventing this requires first writing a derivative of $G_z(\cdot)$ with respect to $z$ as a sum of terms that contain either $G_z$ or a derivative of $G_z$ with respect to the integration variable (which is the argument of $G_z(\cdot)$ in (3.24)).  After doing that, then integrate by parts to rewrite the derivative of the left-most term in (3.24) with respect to $z$ as a sum of integrals whose integrands involves $G_z$ but not its derivatives.  (This last integration by parts step moves the derivative of $G_z$ with respect to its argument off of $G_z$ and onto the product of the functions $\bar\omega_z$, $\ff_0$ and the area 2-form that defines the integration measure.) 	
     
     Note that the preceding issue doesn’t arise with regards to the right-most term on the right-hand side of (3.24) because the argument of $G_z$ and $dG_z$ in the integration is uniformly far from $z$.   
     
     \section{More than 4 points of discontinuity} 	
     
     There exists a set $Z\subset S^2$  with 8 points and a corresponding real line bundle $\mathcal{I} \to S^2-Z$ and harmonic section of $\mathcal{I}$ that is described by (1.1).   The points in this case are the intersections of $S^2$  with the \it lines \rm through the vertices of the inscribed tetrahedron.  The argument for the existence of the data $(Z, \mathcal{I}, \ff)$ in this case is virtually identical to the arguments in the preceding sections.  Alternately:  The points of $Z$ can be viewed as the vertices of a cube inscribed in $S^2$  centered at the origin in $\mathbb{R}^3$.  The section $\ff$ obeying (1.1) is then found using the same arguments as in Sections 2-4 but for the replacement of the group of symmetries of the tetrahedron with the group of symmetries of the cube (the octahedral group). The desired eigensection $\ff$ is the minimizer of the energy functional $\mathcal{E}$ (depicted in (2.2) on an analog of the space $\mathcal{T}_0$ that is defined in this case using the group of orientation preserving symmetries of the cube.  To elaborate:  Let $G$ now denote the octahedral group, the subgroup of $SO(3)$ that preserves the inscribed cube. Let $Z$ denote the set of vertices of the inscribed cube and let $\mathcal{I}$ denote the corresponding real line bundle on $S^2-Z$.  The group $G$’s action on $S^2$ is covered by an isometric action of $\{1,-1\}\times G$ on the bundle $\mathcal{I}$ that covers the action of $G$ on $S^2$, and which is defined so that the element ($-1$, identity) acts on $\mathcal{I}$ as multiplication by the real number $-1$ on the fibers of $\mathcal{I}$.  When $p \in Z$, use $a_p$ to denote the element in $G$ that acts as the clockwise $SO(3)$ rotation by ${2\pi\over 3}$ on the oriented axis along the ray from the origin to $p$.  The set $\{(-1, a_p): p \in Z\}$ generate $\{1, -1\} \times G$ .  Figure 2 schematically depicts the action of the element $(-1, a_1) \in  \{1, -1\} \times G$.  (The figure is explained in detail below.)  Let $\mathcal{T}_0$ denote the space of sections of $\mathcal{I}$ with the property that each group element from the set $\{(-1, a_p): p \in Z\}$ acts as multiplication by $-1$ on the section.  The fact that this space is non-trivial follows from the pattern of $\pm$ signs in Figure 2.  (This is explained momentarily.)   A  minimizing sequence for $\mathcal{E}$ in this current version of $\mathcal{T}_0$ will converge to an element in $\mathcal{T}_0$ which is the desired eigensection section $\ff$.

\begin{figure}[H]
 \centering
 \includegraphics[width=.7\textwidth]{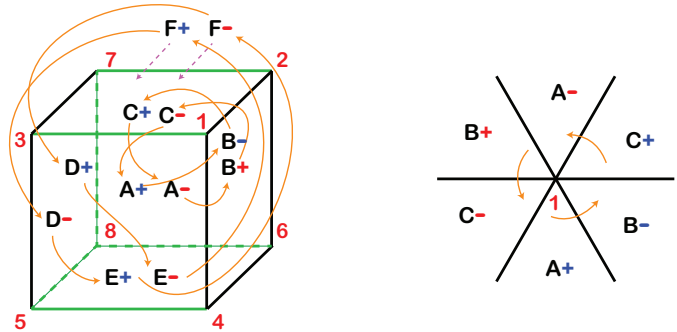}
\caption{The branched cover of the cube and action of $(-1, a_p)$.
}\label{fig:oct}
\vskip -2mm
\end{figure}

	To explain Figure 2:  Let $\mathcal{C}$ denote the cube and let $\mathcal{C}^*$ denote the 2-fold branched cover of the cube, branched over the eight vertices.  The space $\mathcal{C}^*$ is a surface of genus 3.  The complement of the branch points in $\mathcal{C}^*$ is the set of unit length elements in the line bundle $\mathcal{I}$.  The faces of the cube are labeled by capital letters, and their inverse images in $\mathcal{C}^*$ are labeled by a capital letter with an extra $+$ or $-$ label. The left-hand sketch in Figure 2 also
depicts the tiling of $\mathcal{C}^*$ by squares if it is understood that red-colored $\pm$ and blue colored $\pm$ are on different sheets of the cover. (The branch cuts in this depiction of $\mathcal{C}^*$ are indicated by the green edges.)   The vertices of the cube are labeled by the integers in the set $\{1, \ldots, 8\}$.   Each vertex in $\mathcal{C}^*$ has 6 incident squares as indicated in the right-hand drawing of Figure 2.  A rotation by ${2\pi\over 3}$ around a vertex $p$ of the cube is covered by the action of $(-1, a_p)$ on $\mathcal{C}^*$.  The action of $(-1, a_1)$ on a neighborhood of vertex 1 in $\mathcal{C}^*$ is depicted by the right-hand sketch in Figure 1. The orange arrows
in the left-hand sketch indicate how the $\pm$ labeled squares in $\mathcal{C}^*$ are permuted by this
action.  With regards to this action and $\mathcal{T}_0$:  This action moves each square in $\mathcal{C}^*$ to an adjacent square.  The fact that $\mathcal{T}_0$ is non-empty follows from two facts.  The first is that the elements from $\{(-1, a_p): p \in Z\}$ move squares to adjacent squares.  The second is that the inverse images in $\mathcal{C}^*$ of the sides of the cube can be labeled by $\pm$ signs so that no two adjacent squares in $\mathcal{C}^*$ have the same sign label.  To give sign labels with this property, start with the squares in $\mathcal{C}^*$ that are incident on vertex 1.  The right-hand drawing in Figure 2 indicates how to label these squares with $\pm$ labels so that no two adjacent squares have the same label.  (The letter labels of the squares in $\mathcal{C}^*$ are determined by the projection map to the cube. For example, the neighbors of the squares in $\mathcal{C}^*$ incident to vertex 1 and labeled by  $A$  must be labeled by $B$ and $C$.)  Once the signs are set for the squares incident to vertex 1, then there is a unique sign assignment to the remaining squares that makes no two with the same sign adjacent.  Indeed, the signs of the two $F$ labeled squares are a priori determined by examining vertex 2 in $\mathcal{C}^*$ which has the two $F$ labeled squares plus the sign labeled $C\pm$ and $B\pm$.  Likewise, the sign labels of the two $E$ squares in $\mathcal{C}^*$ can be determined from the fact that these $E$ squares are incident to vertex 4 in $\mathcal{C}^*$ where the four other incident squares (the $A\pm$ and $B\pm$ labeled squares) are already labeled by signs.  The signs of the two $D$ labeled squares can be fixed by examining the incident squares to vertex 3 where the other four incident squares ($A\pm$ and $C\pm$) are already sign labeled.  The fact that these sign labels are consistent with the requirement that no two adjacent squares have the same sign can be checked by examining the behavior of the signs around vertices 5, 6 and 7 which each has only one of $A\pm$, $B\pm$ or $C\pm$ incident.  Having checked 5, 6 and 7, then the consistency at vertex 8 follows automatically.
	
	One can also find a set $Z \subset S^2$ with 12 elements, these being the vertices of a regular icosahedron inscribed in the unit radius sphere in $\mathbb{R}^3$.  These vertices are permuted by the icosahedral subgroup of $SO(3)$.  Let $G$ now denote this group.  As explained below, the group $\{1, -1\} \times  G_0$ acts isometrically on the corresponding version of the line bundle $\mathcal{I}$ so as to cover the action of $G$ and so as to have the two crucial properties:  First, the element ($-1$, identity) acts as multiplication by $-1$ on $\mathcal{I}$.  To state second, let $\mathcal{C}$ denote the icosahedron and let $\mathcal{C}^*$ denote the 2-fold branched cover of $S^2$ with branch loci $Z$.  (Keep in mind that the complement in $\mathcal{C}^*$ of the branch loci is the set of unit length elements in $\mathcal{I}$.)  The inverse images in $\mathcal{C}^*$ of the faces of the icosahedron (which are triangles; see Figure 3) can be labeled by $\pm$ signs (this labeling is depicted schematically in Figure 4) so that no two adjacent triangles in $\mathcal{C}^*$ have the same sign label and so that the two inverse images of any face in $\mathcal{C}$ have different signs.  With this understood, suppose for the moment that $p$ is in $Z$ (a vertex of the icosahedron).  Introduce by way of notation $a_p$ to denote the element in $G_0$ that acts as a ${2\pi \over 5}$ clockwise rotation with the axis being the ray from the origin to $p$.  Then the element $(-1, a_p)$ in the group $\{1, -1\} \times  G_0$ acts on $\mathcal{C}^*$ so as to move any given triangle to an adjacent one (this action is depicted schematically in Figure 4).  (This element generates a cyclic subgroup of order 10 whose fifth power is the element ($-1$, identity).)  As in the case of the tetrahedron and the cube, the preceding fact implies that the space $\mathcal{T}_0$ of sections of $\mathcal{I}$ which change sign under the action of any element from $\{(-1, a_p): p \in Z\}$ is non-trivial.  With that understood, then arguments much like those in Proposition 2.2 find a normalized eigensection of the Laplacian in the space $\mathcal{T}_0$ whose norm near any $p\in  Z$ is $\mathcal{O}((|x-p|)^{5\over 2}$) and whose differential has norm $\mathcal{O}((|x-p|)^{3\over2})$ near $p$.
	
There is also a version of $(Z, \mathcal{I}, \ff)$ with $Z$ having 20 points, which are the norm 1 points on the rays from the origin to the midpoints of the faces of the regular icosahedron.  The section $\ff$ in this case is mapped to $-1$ times itself by the generators of the order 6 cyclic groups in $\{1, -1\}\times G$ that map to the order 3 subgroups of the icosahedral group that preserve a given face of the icosahedron.

\begin{figure}[H]
 \centering
 \includegraphics[width=.9\textwidth]{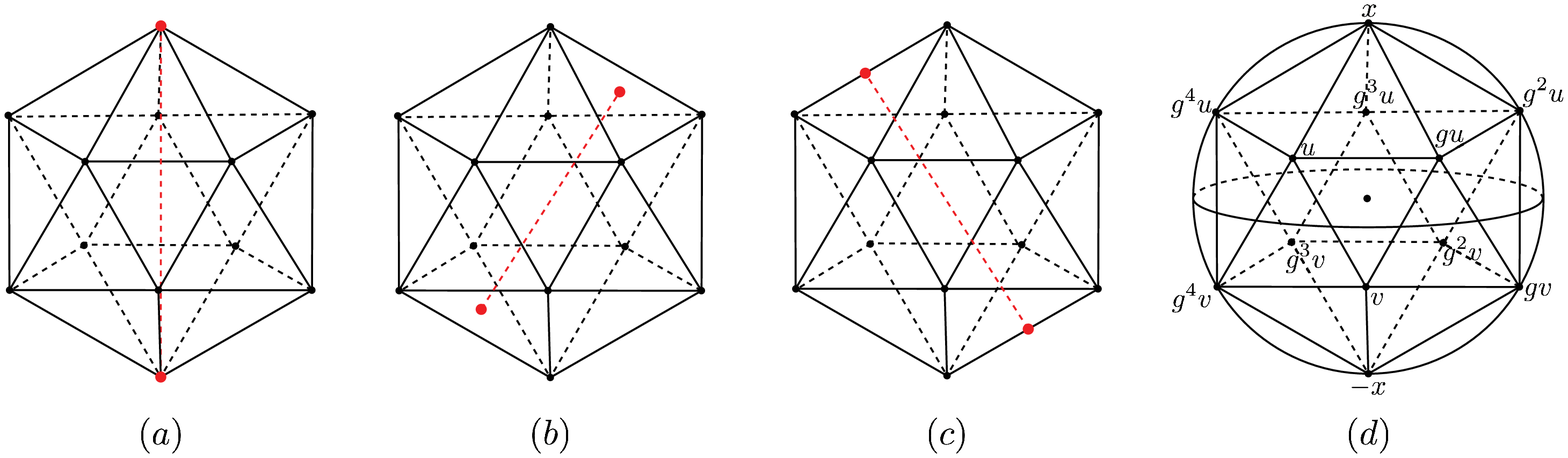}
\vskip 3mm
\caption{Symmetries in an icosahedral group. (a) The $C_5$ axes: symmetries passing through opposite vertices. (b) The $C_3$ axes: symmetries passing through opposite faces. (c) The $C_2$ axes: symmetries passing through middle points of opposite edges.
(d) Three orbits of the symmetry group of a regular icosahedron.
}\label{fig:ico}
\vskip -7mm
\end{figure}

	Figure 3 displays the symmetries in the icosahedral group.  By way of a summary:  This is an order 60 subgroup of the group of $SO(3)$ rotations of $\mathbb{R}^3$ (it is also isomorphic as an abstract group to the alternating group of even permutations of five objects). 
	It appears in $SO(3)$ as the subgroup of rotations that preserve the regular icosahedron.  With this view in mind, some relevant features of the icosahedral  group are depicted by Figure \ref{fig:ico}.  To explain the figure, keep in mind first that a non-trivial rotation of $\mathbb{R}^3$ fixes precisely 2 points, one being the antipodal of the other.  It then acts as a rotation of the plane perpendicular to the line through these points (this line is called the \it axis \rm of the rotation.)  An axis of rotation for an element in the icosahedral group is denoted by $C_n$ where $n \in \{2, 3, 5\}$ is the order of the group element.  (This is the smallest positive integer $n$ such that rotation about the relevant axis by ${2\pi \over n}$ gives an equivalent configuration.)  The various $C_n$'s for the icosahedral group are as follows: 

\begin{enumerate}
\item There are six $C_5$ axis; these are the axes through antipodal vertices.  One is 
       depicted in Figure \ref{fig:ico}a.
\item  There are ten $C_3$ axis; these are the axes through antipodal faces.  One is 
       depicted in Figure \ref{fig:ico}b.
\item There are fifteen $C_2$ axis; these are the axes that bisect antipodal edges.  One is 
       depicted in Figure \ref{fig:ico}c.
\end{enumerate}
By way of an example, Figure \ref{fig:ico}d depicts the orbit of the vertices under the $C_5$ rotation that fixes the vertices labeled $x$ and $-x$.  The orbits are as follows:
$$\{x, -x, u, gu, g^2u, g^3u, g^4u, v, gv, g^2v, g^3v, g^4v\}.$$

Let $G$ again denote the icosahedral group.  Two key points were noted above about the group $\{1, -1\} \times G$ and its action on $\mathcal{C}^*$ (the 2-fold branched cover of $S^2$ with branch loci $Z$, the vertices of the icosahedron).  These both concern the labeling of the inverse images in $\mathcal{C}^*$ of the faces of the tetrahedron:  That these faces can be labeled by $\pm$ signs so that no two adjacent triangles have the same sign label and so that the two inverse images of any given face of the icosahedron have different sign labels.  Such a labeling is illustrated schematically by Figure 4:

\begin{figure}[H]
 \centering
  \includegraphics[width=\textwidth]{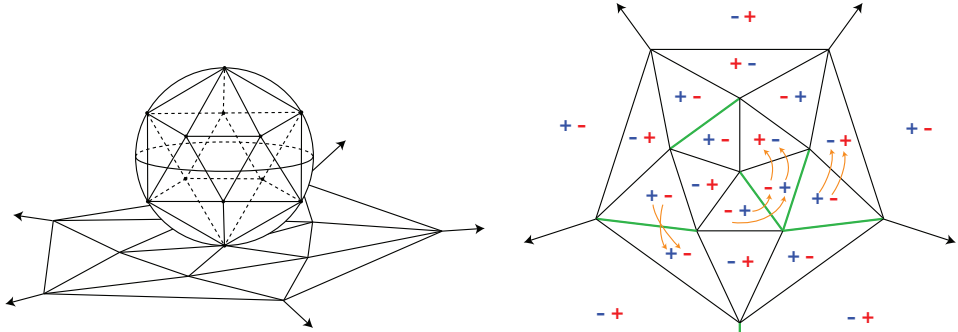}
\caption{The labeling of the triangles in $\mathcal{C}^*$}\label{fig:proj}
\vskip -3mm
\end{figure}

To explain Figure 4:  The left-hand diagram in Figure 4 shows a stereographic projection of the icosahedron to $\mathbb{R}^2$ whose image is depicted in the right-hand diagram.  Each face in the right-hand diagram has two labels, $+$ and $-$.  These are the sign labels of the inverse image of the face in $\mathcal{C}^*$.  Moving in $\mathcal{C}^*$ to an adjacent triangle can be viewed as moving the corresponding triangle in the right-hand diagram to an adjacent triangle with the proviso that the color of the sign (blue or red) must stay the same except when crossing a green edge where the color must change.  (The color of the $\pm$ signs in the right-most sketch of Figure 4 distinguish the two sheets of the 2-fold cover. The green edges in the right-hand sketch of Figure 4 signify a placement of the branch cuts on the icosahedron.) Keeping this rule in mind, then the diagram exhibits a $\pm$ labeling of the triangles in $\mathcal{C}^*$ with the desired property. (One can also prove that such a labeling exists by starting with a labeling of the ten triangles in $\mathcal{C}^*$ adjacent to a given vertex and then moving from vertex to vertex much as was done in the case of the cube.)   The orange arrows in the right-hand diagram indicate how some of the sign labeled triangles in $\mathcal{C}^*$ move under the action of the element $(-1, a_p)$ when $p$ is the central vertex in the right-hand diagram.

\vskip 15mm

\begin{appendices}
\section*{Appendix}

	The purpose of this appendix is to state and prove the following proposition:\\

\noindent\bf Proposition A: \it
 Let $(Z, \mathcal{I}, \ff)$ denote a data set with $Z$ being two distinct points in $S^2$, with $\mathcal{I} \to S^2-Z$ denoting the non-trivial real line bundle, and with $\ff$ denoting a Laplace eigensection of $\mathcal{I}$ that obeys the conditions in (2.1).  Then the points in $Z$ are antipodal and $\ff$ is the real part of $az^{k+{1\over2}}$, where $a$ is a non-zero complex number, and $z$ is a complex Euclidean coordinate on the plane perpendicular to the line between the two points of $Z$.\\ 

\noindent\textbf{Proof of Proposition A:} \rm
If the two points in $Z$ are antipodal, then the $S^1$ action on $S^2$ that rotates the sphere around the line in $\mathbb{R}^3$ through the two points is covered by an isometric fiber preserving action of $S^1$ on $\mathcal{I}$.  (With $S^1$ viewed as $\mathbb{R}/2\pi \mathbb{Z}$, the action of $\pi\in\mathbb{R}/2\pi \mathbb{Z}$ on $\mathcal{I}$  covers the $2\pi$ rotation of $S^2$, and it acts as multiplication by $-1$ on the fibers.)  Because this $S^1$ action is isometric, the eigensections for the Laplace operator can be found using a standard Fourier series separation of variables.  Doing this leads to the form of the eigensections that is described in the proposition. The proof that the points of $Z$ must be antipodal has five parts.  There is also a Part 6 with an extra parenthetical remark.\\

\it	Part 1: \rm  Fix an even number of distinct points in $S^2$ for $Z$ and then construct the associated real line bundle $\mathcal{I}$ on  $S^2 - Z$.   Let $\fp$ and $\fq$ denote Laplace eigensections of $\mathcal{I}$ with the same eigenvalue.  With regards to $\fp$:  Assume that $\fp$ and its derivative $d\fp$ can be written near each $p \in Z$ as the respective real parts of 

$$a_{p} \sqrt{z} + \mathcal{O}(|z|^{3\over2})   \quad    and   \quad  {1\over 2}  a_{p}  {1\over \sqrt{z}}  dz + \mathcal{O}(|z|^{1\over2})$$ 
\hfill (A.1)\\
with $a_p$ being a complex number.  With regards to $\fq$:  Assume that $\fq$ and its derivative $d\fq$ can be written near each $p \in Z$ as the respective real parts of 
$$c_{p}  {1\over \sqrt{z}} dz + \mathcal{O}(|z|^{1\over2})  \quad    and \quad   -{1\over 2} c_{p} {1\over z^{3\over2}} dz + \mathcal{O}(|z|^{-{1\over2}}).$$
\hfill (A.2)\\

 	Let $*$ denote the round sphere's Hodge star operator.  The 1-form 

$$* (\fq d\fp -  \fp d\fq) $$
\hfill (A.3)\\
is necessarily a closed 1-form.  This is because $d* d\fp$ is $-\mathcal{E} \fp$ times the area 2-form and, likewise, $d* d\fq$ is $-\mathcal{E} \fq$ times this same 2-form.  (And because $d\fg \wedge * d\fk = d\fk \wedge * d\fg$ when $\fk$ and $\fg$ are any two sections of $\mathcal{I}$.)  With the preceding understood, fix a small positive number to be denoted by $\rho$ and reintroduce the function $\chi_{\rho}$ from Part 1 of the proof of Proposition 2.1 in Section 3.  Then
$$\int_{S^2} \chi_\rho d*(\fq d\fp - \fp d\fq) = 0.$$
\hfill (A.4)	\\
This last identity leads to a bilinear relation between the various $p \in Z$ of the complex numbers $a_p$ and $c_p$ that appears in (A.1) and (A.2).  To obtain the desired relation, integrate by parts in (A.4) to write the left-hand side integral as
$$-\int_{S^2} d\chi_\rho \wedge * (\fq d\fp - \fp d\fq)$$ 
\hfill (A.5)\\
This is a sum of integrals indexed by the points in $Z$ with any given $p \in Z$ contribution to the sum being an integral whose integrand is supported where the distance to $p$ is less than 2$\rho$.  As a consequence, (A.4) and (A.5) can be used to evaluate the contribution from each $p \in Z$ up to leading order in $\rho$.  The result of doing so is a sum whose $\rho \to 0$ limit is
$$\sum_{p\in Z} (a_pc_p + \bar a_p \bar c_p ) = 0.$$
\hfill (A.6)\\
(The two contributions $\fq d \fp$ and $-\fp d\fq$ give the same contribution to (A.6) because the two appearances of $a_p$ in (A.1) have the same sign whereas the two appearances of $c_p$ in (A.2) have opposite signs.)\\

\it Part 2: \rm  With (A.6) in mind, this part of the proof explains how any given Laplace eigensection generates a set of Laplace eigensections with each having the same eigenvalue as the given one.   To this end, let $L_1$, $L_2$ and $L_3$ denote the generators of the $SO(3)$ rotations about the respective $x_1$, $x_2$ and $x_3$ axis.  Thus,
\begin{align*}
L_1 = x_2{\partial \over \partial x_3}  - x_3 {\partial \over \partial x_2} \quad    and  \quad  
L_2 = x_3{\partial \over \partial x_1}  - x_1 {\partial \over \partial x_3} \quad    and  \quad 
L_3 = x_1{\partial \over \partial x_2}  - x_2& {\partial \over \partial x_1}.\\
&\text{(A.7)}
\end{align*}

These operators can be viewed as acting on the space of functions on $S^2$.  Moreover, supposing that $Z$ is a set of some even number of distinct points in $S^2$ and $\mathcal{I} \to S^2 - Z$  is the corresponding real line bundle, these operators also act on the space of sections of $\mathcal{I}$.  And, in any of these incarnations, the operators in (A.7) commute with the spherical Laplacian (which can be written as $L_1L_1 + L_2L_2 + L_3L_3$.)
 
The preceding facts have the following implication:  If $\ff$ is a Laplace eigensection of $\mathcal{I}$ with eigenvalue $\mathcal{E}$, then so is any constant linear combination from the set 

$$\{L_a\ff\}_{a=1,2,3}\cup \{L_aL_b\ff\}_{a,b=1,2,3} \cup \cdots$$
\hfill (A.8)\\
where the elements in the unwritten part have the form $L_{a_1}\cdots L_{a_m}\ff$  with $m \geq 3$.\\

\it Part 3:  \rm Fix $p \in Z$ (with $Z$ as just described) and let $z$ denote a holomorphic coordinate for $S^2$ that is defined near $p$ with norm $|dz|$ at $p$ equal to $\sqrt{2}$.  Suppose that $\ff$ is an eigensection of $\mathcal{I}$ for the Laplacian with the integrals of $|\ff|^2$ and $|d\ff|^2$ being finite.  The arguments from Parts 2-5 of the proof of Proposition 2.2 can be used to see that $\ff$ and $d\ff$ near $p$ can be written as the real parts of 
$$ a z^{k+{1\over2}} + \mathcal{O}(|z|^{k+{3\over2}}) \quad  and  \quad  {1\over 2} a z^{k-{1\over2}} dz + \mathcal{O}(|z|^{k+{1\over2}}) .$$
\hfill (A.9)\\
with $k$ being a non-negative integer and with $a \in \mathbb{C}-\{0\}$.  Moreover, the Laplace eigensection $L_a\ff$ and its exterior derivative near $p$ are the respective real parts of 
$$\left(k+ {1\over 2}\right) a z^{k-{1\over 2}} (L_az)_{z=0} + \mathcal{O}(|z|^{k+{1\over 2}}) \ and \ \left(k^2  - {1\over 4} \right) a z^{k-{3\over 2}} (L_az)_{z=0}dz + \mathcal{O}(|z|^{k-{1\over 2}})$$
\hfill (A.10)\\

	As a consequence of these asymptotics, what is said in Part 1 can be invoked using the eigensection $\ff$ for $\fp$ and $L_a\ff$ for $\fq$ (with any choice of $a \in \{1, 2, 3\}$).  In this case, the complex number $a_p$ is zero if the integer $k$ that appears in $p$'s version of (A.9) is positive; and it is the complex number $a$ in (A.9) if $k = 0$.  Meanwhile, $c_p$ is also zero if $k > 0$ and it is equal to $\left(k+{1\over 2} \right) a_p (L_zz)_{z=0}$ if $k = 0$.  With this understood, (A.6) asserts the vanishing of the real part of 
$$\sum_{p\in Z} a_p^2(L_zz)_{z=0}.$$
\hfill (A.11)\\
	In general, if $m$ is a positive integer, and if all $ p \in Z$ versions of (A.9)’s integer $k$ are no smaller than $m$, then what is said in Part 1 can be invoked using $\fp = L_{a_1}\cdots L_{a_{m}}\ff$ for any choice of $(a_1, \ldots, a_m) \in \{1, 2, 3\}^m$ and using $\fq = L_{a_{m+1}} \cdots L_{a_{2m+1}} \ff$  for any choice of $(a_{m+1}, \ldots , a_{2m+1}) \in \{1, 2, 3\}^{m+1}$.  In this event, $a_p$ is zero if $p$’s version of (A.9) has $k > m$; and up to a $p$-independent positive factor, $a_p$ is equal to
 $a \prod_{j=1}^m (L_{a_j}z)_{z=0}$  otherwise.  By the same token, $c_p$ is zero if $k > m$; and up to a $p$-independent, positive factor, $c_p$ is equal to $a \prod_{j=m+1}^{2m+1} (L_{a_j}z)_{z=0}$ when $k = m$.  The relation in (A.6) in this case asserts the vanishing of the real part of the complex sum
$$\sum_{p\in Z}a_p^2 \prod_{j=1}^{2m+1} (L_{a_j}z)_{z=0}.$$
\hfill (A.12)\\
(Note in particular that the larger the value of $m$, the greater the number of relations that must be satisfied by the squares of the leading order coefficients of $\ff$ near the points in $Z$.)\\
 
 \vskip 2mm

\it Part 4: \rm To say more about (A.12), it proves useful to make a coherent choice for the local coordinate $z$ near each point in $Z$.  To this end, no generality is lost by assuming that the south pole of $S^2$ is not a point in $Z$.  Stereographic projection from the south pole now identifies the complement of that point with $\mathbb{C}$ as follows:  The complex Euclidean coordinate on $\mathbb{C}$ is denoted by $u$, and it is given in terms of the Euclidean coordinates $(x_1, x_2, x_3)$ on the sphere by the rule whereby
$$u = {x_i + i x_2 \over 1 + x_3}.$$
\hfill (A.13)\\

If $p \in Z$, let $u_p$ denote the value of $u$ at $p$.  The complex coordinate 
$$z ={2\over (1 + |u_p|^2) } (u - u_p)$$ 
\hfill (A.14)\\
is zero at $p$ and the norm of $dz$ at $p$ is $\sqrt{2}$.  It can therefore be used in (A.9) and (A.10).  With regards to $L_az$:
\begin{itemize}
\item $L_1z =  - {i\over (1 + |u_p|^2)}(1-u_p^2)$.
\item $L_2z =  {1\over (1 + |u_p|^2)}(1+u_p^2)$.
\item $L_3z =  {2i u_p\over (1 + |u_p|^2)}$.
\end{itemize}
\hfill (A.15)\\
Note that $(L_1z)^2 + (L_2z)^2 + (L_3z)^2 = 0$.  This implies that the constraints that are implied by the vanishing of the real parts of the expressions in (A.12) for a given $m \geq 1$ are not all linearly independent.  For example, there are only 7 real-valued constraints when $m = 1$.  \\

\vskip 2mm

\it	Part 5: \rm  Now suppose that $Z$ has two points.  One can be moved to the $u = 0$ point by an $SO(3)$ rotation of $S^2$ and, if they are not antipodal, the other can be moved to a point where $u$ is real and positive.  Suppose that the respective versions of $k$ that appear in (A.9) for these points are positive.  Denote the version of $a$ for the $u = 0$ point as $a_0$ and the version of $a$ for the $u > 0$ point as $a_1$.  

The expressions in (A.12) with $m = 1$ involve three versions of $L_az$.  Those with $a = 3$ involve only the $u > 0$ point because $L_3z = 0$ when $u_p = 0$.  Those lead to the following conditions:  Taking $a_1 = a_2 = a_3 = 3$ leads to the constraint
$$a_1^2 - \bar a_1^2  = 0,$$
\hfill (A.16)\\
which is to say that $a_1^2$ is a real number.  The constraint with $a_1 = a_2 = 3$ and $a_3 = 1$ leads again to (A.16) unless $u_p = 1$ (in which case the constraint is vacuous), and the constraint with $a_1 = a_2 = 3$ and leads to  respective constraints
$$a_1^2 + \bar a_1^2 = 0,$$
\hfill (A.17)\\
which is to say that $a_1^2$ is real.  As a consequence of both of these, $a_1$ must vanish.  Then, the respective constraints with $a_1 = a_2 = a_3 = 1$ and $a_1 = a_2 = a_3 = 2$ require that $a_0^2$ to be first real and then imaginary, so it too must vanish.  

	It follows as a consequence that the integer $k$ at both points must be greater than 1.  Similar arguments (using induction) show that $k$ can not be any positive integer, which is nonsensical because it runs afoul of what is said in Part 3.\\
	
	\it Part 6:  \rm By way of a parenthetical remark:  In the case of the tetrahedron, it is an exercise to check that the 7 real-valued constraints can be satisfied (and likewise with the other cases from Section 5).  In the tetrahedral case, the four points in $Z$ are listed in (2.4).  If the left-most point (1, 0, 0) is the $u = 0$ point in $\mathbb{C}$, then the constraints are obeyed if all four of the $a_p$’s have the same norm with the $u = 0$ version obeying $a_p^2 < 0$ and the other three points in $Z$ obeying $a_p^2 > 0$. (The fact that the $u = 0$ vertex of the tetrahedron has $a_p^2 < 0$ whereas the other vertices have $a_p^2 > 0$ is consistent with regards to a given eigensection being equivariant with respect to the action on $\mathcal{I}$ of the product of $\{1, -1\}$ with the tetrahedral group.  This sign change is due to the behavior of the differential of the stereographic projection map.)

\end{appendices}

\bibliographystyle{plain}
\bibliography{bib.bib}

\vskip 10mm

\begin{tabular}{ ll }
C. H. Taubes: &Department of Mathematics\\
& Harvard University\\
& Cambridge, MA, 02138\\
 &\url{chtaubes@math.harvard.edu}\\
 &\\
 Y. Wu: & Center of Mathematical Sciences and Applications\\
 & Harvard University\\ 
& Cambridge, MA, 02138\\
&	\url{ywu@cmsa.fas.harvard.edu}\\
\end{tabular}
\\
 
\end{document}